\documentclass[10pt,journal]{IEEEtran}

\usepackage{array}

\usepackage{graphicx}
\usepackage[usenames]{color}

\usepackage{amsmath}  
\usepackage{amssymb}  

\usepackage{theorem}  
\usepackage{cite}     
\usepackage{comment}  

\usepackage{url}

\usepackage{amscd}
\usepackage{latexsym}

\usepackage{pict2e}

\title{Two-Dimensional Patterns with Distinct Differences -- Constructions, Bounds, and Maximal Anticodes }

\author{Simon~R.~Blackburn, Tuvi~Etzion, Keith~M.~Martin and
Maura~B.~Paterson%
\thanks{This work was supported in part by the
Israel Science Foundation (ISF) under grant no.\ 230/08, as well as by
EPSRC grants EP/D053285/1 and EP/F056486/1.}%
\thanks{S.R.~Blackburn,
K.M.~Martin and M.B.~Paterson are with the Department of Mathematics,
Royal Holloway, University of London, Egham, Surrey TW20 0EX.
T.~Etzion is with the Computer Science Department, Technion--Israel
Institute of Technology, Haifa 32000, Israel.}}

\newcommand{\cA}{{\cal A}}
\newcommand{\cB}{{\cal B}}

\newcommand{\cS}{{\cal S}}

\newcommand{\cX}{{\cal X}}

\newcommand{\bZ}{\mathbb{Z}}

\newcommand{\DD}{\mathrm{DD}}
\newcommand{\oDD}{\overline{\mathrm{DD}}}

\newcommand{\highsup}[1]{\raisebox{0.35ex}{\kern 1pt $\scriptstyle {#1} $}}

\makeatletter
\DeclareRobustCommand{\sbinom}{\genfrac[]\z@{}}
\makeatother
\newcommand{\G}[2]{\sbinom{{#1}\kern-1pt}{{#2}\kern-1pt}}
\newcommand{\Gq}[2]{\sbinom{{#1}\kern-1pt}{{#2}\kern-0.5pt}}

\newcommand{\Z}{{\mathbb Z}}
\newcommand{\R}{{\mathbb R}}
\newcommand{\Q}{{\mathbb Q}}

\newcommand{\be}[1]{\begin{equation}\label{#1}}
\newcommand{\ee}{\end{equation}}

\newcommand{\bc}{\begin{center}}
\newcommand{\ec}{\end{center}}

\DeclareMathAlphabet{\mathbfsl}{OT1}{cmr}{bx}{it}
\newcommand{\uuu}{\kern-1pt\mathbfsl{u}\kern-0.5pt}
\newcommand{\vvv}{\kern-1pt\mathbfsl{v}\kern-0.5pt}

\newcommand{\sP}{\script{P}}

\newcommand{\sG}{\script{G}}

\newcommand{\Ps}{\smash{{\sP\kern-2.0pt}_q\kern-0.5pt(n)}}
\newcommand{\sPs}{\smash{{\sP\kern-1.5pt}_q(n)}}
\newcommand{\Gr}{\smash{{\sG\kern-1.5pt}_q\kern-0.5pt(n,k)}}
\newcommand{\Grmk}{\smash{{\sG\kern-1.5pt}_q\kern-0.5pt(n,n-k)}}
\newcommand{\Grdk}{\smash{{\sG\kern-1.5pt}_q\kern-0.5pt(2k,k)}}
\newcommand{\Grekappa}{\smash{{\sG\kern-1.5pt}_q\kern-0.5pt(n,e+1-\kappa)}}
\newcommand{\Grtwoekappa}{\smash{{\sG\kern-1.5pt}_q\kern-0.5pt(n,2e+1-\kappa)}}
\newcommand{\Gremkappa}{\smash{{\sG\kern-1.5pt}_q\kern-0.5pt(n,e-\kappa)}}
\newcommand{\Ptwo}{\smash{{\sP\kern-2.0pt}_2\kern-0.5pt(n)}}
\newcommand{\Ptwom}{\smash{{\sP\kern-2.0pt}_2\kern-0.5pt(m)}}
\newcommand{\Ptwonm}{\smash{{\sP\kern-2.0pt}_2\kern-0.5pt(n+m)}}
\newcommand{\Ptwoa}{\smash{{\sP\kern-2.0pt}_2\kern-0.5pt(1)}}
\newcommand{\Ptwob}{\smash{{\sP\kern-2.0pt}_2\kern-0.5pt(2)}}
\newcommand{\Ptwoc}{\smash{{\sP\kern-2.0pt}_2\kern-0.5pt(3)}}
\newcommand{\Ptwod}{\smash{{\sP\kern-2.0pt}_2\kern-0.5pt(4)}}
\newcommand{\Ptwoe}{\smash{{\sP\kern-2.0pt}_2\kern-0.5pt(5)}}
\newcommand{\Ptwokm}{\smash{{\sP\kern-2.0pt}_2\kern-0.5pt(2k-1)}}

\newcommand{\Gn}{\smash{{\sG\kern-1.5pt}_2\kern-0.5pt(n,n{-}1)}}
\newcommand{\Gnq}{\smash{{\sG\kern-1.5pt}_q\kern-0.5pt(n,n{-}1)}}
\newcommand{\Gone}{\smash{{\sG\kern-1.5pt}_2\kern-0.5pt(n,1)}}
\newcommand{\GTwo}{\smash{{\sG\kern-1.5pt}_2\kern-0.5pt(n,k)}}
\newcommand{\Gnk}{\smash{{\sG\kern-1.5pt}_2\kern-0.5pt(n,n{-}k)}}
\newcommand{\Pone}{\smash{{\sP\kern-2.5pt}_2\kern-0.5pt(n{-}1)}}
\newcommand{\Greone}{\smash{{\sG\kern-1.5pt}_q\kern-0.5pt(n,e{+}1)}}
\newcommand{\Gretwo}{\smash{{\sG\kern-1.5pt}_q\kern-0.5pt(n,e{+}2)}}

\newcommand{\myboxplus}{\kern1pt\mbox{\small$\boxplus$}}

\renewcommand{\leq}{\leqslant}

\renewcommand{\geq}{\geqslant}

\newcommand{\script}[1]{{\mathscr #1}}

\newcommand{\Cref}[1]{Co\-rol\-la\-ry\,\ref{#1}}

\theoremstyle{plain}
\theorembodyfont{\normalfont\slshape}

\newtheorem{thm}{Theorem$\!$}
\newenvironment{theorem}{\begin{thm}\hspace*{-1ex}{\bf.}}{\end{thm}}

\newtheorem{prop}[thm]{Proposition$\!$}

\newtheorem{lem}[thm]{Lemma$\!$}
\newenvironment{lemma}{\begin{lem}\hspace*{-1ex}{\bf.}}{\end{lem}}

\newtheorem{cor}[thm]{Corollary$\!$}
\newenvironment{corollary}{\begin{cor}\hspace*{-1ex}{\bf.}}{\end{cor}}

\newtheorem{defi}{Definition$\!$}
\newenvironment{definition}{\begin{defi}\hspace*{-1ex}{\bf.}}{\end{defi}}

\theorembodyfont{\normalfont}

\newtheorem{exam}{Example$\!$}

\begin{document}

\maketitle

\begin{abstract}
A two-dimensional grid with dots is called a \emph{configuration with
distinct differences} if any two lines which connect two dots are
distinct either in their length or in their slope. These
configurations are known to have many applications such as radar,
sonar, physical alignment, and time-position synchronization.  Rather
than restricting dots to lie in a square or rectangle, as previously
studied, we restrict the maximum distance between dots of the
configuration; the motivation for this is a new application of such
configurations to key distribution in wireless sensor networks. We
consider configurations in the hexagonal grid as well as in the
traditional square grid, with distances measured both in the Euclidean
metric, and in the Manhattan or hexagonal metrics.

We note that these configurations are confined inside maximal anticodes
in the corresponding grid. We classify maximal anticodes for each
diameter in each grid. We present upper bounds on the number of dots
in a pattern with distinct differences contained in these maximal
anticodes. Our bounds settle (in the negative) a question of Golomb
and Taylor on the existence of honeycomb arrays of arbitrarily large
size. We present constructions and lower bounds on the number of dots
in configurations with distinct differences contained in various
two-dimensional shapes (such as anticodes) by considering periodic
configurations with distinct differences in the square grid.
\end{abstract}


\section{Introduction}


\PARstart{A}{} {\it Golomb ruler} (or {\it ruler} for short) of order $m$ (also
known as a {\it Sidon set}) is a set $S$ of integers with $|S|=m$
having the property that all differences $a-b$ (for $a,b \in S$, with
$a \neq b$) are distinct. They were first used by Babcock, in
connection with radio interference~\cite{Bab}. The {\it length} of a
Golomb ruler $S$ is the largest difference between any two elements of
$S$. It is easy to show that a ruler of order $m$ has length at least
$\binom{m}{2}$; a ruler meeting this bound is called {\it
perfect}. Golomb has shown that no perfect ruler exists with order
greater than four~\cite{Gol72}.  The problem of finding the shortest
possible length of a Golomb ruler of a given order has been widely
studied; no general solution is known, but optimal rulers have been
determined for orders less than 24 (see \cite{Shearer} for
details). The elements of a Golomb ruler can be taken to represent
marks (`dots') on a physical ruler occurring at integer differences
from each other. The fact that the differences are all distinct
implies that if a Golomb ruler is placed on top of a second,
identical, ruler then at most one mark from the upper ruler will
coincide with a mark from the lower ruler, unless they are exactly
superimposed. Golomb rulers arise in the literature from both
theoretical and practical aspects
(see~\cite{ASU,Bab,LaSa88,GoTa82}). It is well known that the largest
order of a ruler of length $n$ is $\sqrt{n}+o(\sqrt{n})$,
see~\cite{ASU,LaSa88,RoBe67}.

There are various generalizations of one-dimensional rulers into
two-dimensional arrays. One of the most general was given by
Robinson~\cite{Rob85}. A two-dimensional ruler is an $n \times k$
array with $m$ dots such that all $\binom{m}{2}$ lines connecting two
dots in the array are distinct as vectors, i.e., any two have either
different length or slope. These arrays were also considered
in~\cite{Rob97,Rob00}. The case where $n=k$ was first considered suggested by Costas and investigated by
Golomb and Taylor~\cite{GoTa82}. Costas considered the case when $n=k$
and each row and each column in the array has exactly one
dot~\cite{GoTa82}. These arrays have application to a sonar
problem~\cite{Cos75,GoTa82} and also to radar, synchronization, and
alignment; they are known as Costas arrays. Sonar sequences
are another class of arrays mentioned in~\cite{GoTa82}, where $m=k$ and
each column has exactly one dot; see~\cite{GRT87,Gam87,MGT93}.

Other two-dimensional generalizations of a Golomb ruler have been
considered in the literature, but do not have direct connection to
our current work. For the sake of completeness we will mention them.
A {\it radar array} is an $n \times k$ array with exactly one dot
per column such that there are no two lines connecting two disjoint
pairs of dots, occupying the same rows, which have the same length
and slope. Radar arrays were defined in~\cite{GoTa82} and considered
in~\cite{BlTi88,HaZe97,Rob85,ZhTu94}. Arrays in which all lines have
distinct slopes were considered in~\cite{EGRT92,PeTa00,Zha93}.
Arrays in which the Euclidean distances of any pair of lines are
distinct were considered in~\cite{LeTh95}.


\IEEEpubidadjcol

The examples above are concerned with dots in an (infinite) square
grid that are restricted to lie in a given line segment, square or
rectangle. More generally, we can define a set of dots in a grid to be
a \emph{distinct difference configuration} (DDC) if the lines
connecting pairs of dots are different either in length or in
slope. Having surveyed the known structures of two-dimensional
patterns with distinct differences, it seems that the following
natural question has not been investigated: what is the maximum number
of dots that can be placed on a two-dimensional square grid such that
all lines connecting two dots are different either in their length or
their slope and the distance between any two dots is at most $r$? In
other words, rather than considering the traditional rectangular
regions of the square grid, we consider dots which lie in maximal
anticodes of diameter $r$. (An \emph{anticode of diameter $r$} is a
set of positions in the grid such that any pair of positions are at
distance at most~$r$.  See Section~\ref{sec:anticodes} for details.)
We will consider two notions of distance in the square grid: Manhattan
and Euclidean. We also consider distinct difference configurations in
the hexagonal grid, using hexagonal distance (also known as Manhattan
distance) and Euclidean distance.

Our motivation for considering these configurations comes from a new
application to key predistribution for wireless sensor networks. We
considered in~\cite{BEMP} a key predistribution scheme based on DDCs
in general, and Costas arrays in particular. A DDC $\cA$ with $m$ dots
was shifted over the two-dimensional square grid. For each shift we
assigned the same key to the $m$ entries of the two-dimensional grid
which coincide with the $m$ dots of $\cA$.  In \cite{BEMP} we noted
that a Costas array is a DDC, and gave examples of DDCs with small
numbers of dots.  However, the questions of finding more general
constructions, and providing bounds on the number of dots in such
configurations, were left open; it is these issues that are addressed
by the results of this paper.  Other properties of DDCs motivated by
this application are considered in~\cite{BEMP1}.

The rest of this paper is organized as follows. In
Section~\ref{sec:models} we describe the models on which we will
consider our two-dimensional patterns with distinct differences. We
consider two two-dimensional grids: the square grid and the hexagonal
grid. In the square grid we consider the Manhattan distance and the
Euclidean distance, while in the hexagonal grid we consider the
hexagonal distance and the Euclidean distance. We define the classes
of DDCs we will study, and list optimal examples for small parameter
sizes. In Section~\ref{sec:anticodes} we explain the relation between
DDCs and maximal anticodes. We classify the maximal anticodes when we
use Manhattan distance and hexagonal distance. We also briefly review
some properties of anticodes in $\R^2$ using Euclidean distance: these
properties will allow us to bound the size of an anticode in either
grid when we use Euclidean distance.  In
Section~\ref{sec:upper_bounds} we present upper bounds on the number
of dots in a DDC when we restrict the dots to lie in some simple
regions (`shapes') in the grid. The most important shapes we consider
are the anticodes, in particular the Lee sphere and the hexagonal
sphere. As a consequence of our upper bound, we settle an old question
of Golomb and Taylor~\cite{GoTa84} (on the existence of honeycomb
arrays of arbitrarily large size) in the negative. In
Sections~\ref{sec:TwoDPatterns} and~\ref{sec:lower} we turn our
attention to constructions and lower bounds for the number of dots in
the classes of DDC defined in Section~\ref{sec:models}. We generalize
a folding technique that was used by Robinson~\cite{Rob97} to construct
Golomb rectangles, and provide more good examples by constructing
periodic infinite arrays that are locally DDCs. Our constructions are
asymptotically optimal in the case of the square grid and Manhattan
distance.
\section{Grids, Distances, and DDCs}
\label{sec:models}

We first define some new classes of two-dimensional distinct
difference configurations. We believe that the definitions are very
natural and are of theoretical interest, independently of the
application we have in mind. We will consider the square grid and the
hexagonal grid as our surface. We start with a short definition of the
two models. Before the formal definition we emphasize that we define a
point $(i,j)$ to be the point in column $i$ and row $j$ of either a
coordinate system or a DDC. Hence, rows are indexed from bottom to top
in increasing order; columns are indexed from left to right in
increasing order. (So this is the usual convention for a Cartesian
co-ordinate system, but is not the standard way of indexing the
entries of a finite array.)

\subsection{The two models}

The first model is called the {\it square model}. In this model, a
point $(i,j)\in\Z^2$ has the following four neighbors when we consider the
model as a connected graph:
$$\{(i-1,j),(i,j-1),(i,j+1),(i+1,j)\}.$$ We can think of the points in
$\Z^2$ as being the centres of a tiling of the plane by unit squares,
with two centres being adjacent exactly when their squares share an
edge.  The distance $d((i_1,j_1),(i_2,j_2))$ between two points
$(i_1,j_1)$ and $(i_2,j_2)$ in this model is the Manhattan distance
defined by
$$d((i_1,j_1),(i_2,j_2))=|i_2-i_1|+|j_2-j_1|.$$

The second model is called the {\it hexagonal model}. Instead of the
square grid, we define the following graph. We start by tiling the
plane $\R^2$ with regular hexagons whose sides have length $1/\sqrt{3}$
(so that the centres of hexagons that share an edge are at distance
$1$). The vertices of the graph are the centre points of the
hexagons. We connect two vertices if and only if their respective
hexagons share an edge. This way, each vertex has exactly six
neighbouring vertices.

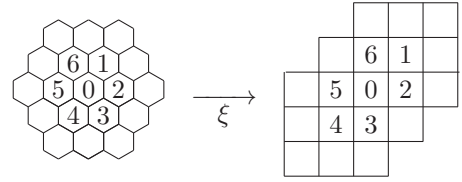
\begin{figure}[t]
\centering
\setlength{\unitlength}{.4mm}
\begin{picture}(150,70)(-10,-40)
\multiput(10,5)(10,0){3}{
\put(-10,0){\line(0,1){5.77350269}}
\put(0,0){\line(0,1){5.77350269}}
\put(0,5.77350269){\line(-1732,1000){5}}
\put(-10,5.77350269){\line(1732,1000){5}}
}
\multiput(5,-3.66033872)(10,0){4}{
\put(-10,0){\line(0,1){5.77350269}}
\put(0,0){\line(0,1){5.77350269}}
\put(0,5.77350269){\line(-1732,1000){5}}
\put(-10,5.77350269){\line(1732,1000){5}}
}
\multiput(0,-12.3206774)(10,0){5}{
\put(-10,0){\line(0,1){5.77350269}}
\put(0,0){\line(0,1){5.77350269}}
\put(0,5.77350269){\line(-1732,1000){5}}
\put(-10,5.77350269){\line(1732,1000){5}}
}
\multiput(5, -20.9810162)(10,0){4}{
\put(-10,0){\line(0,1){5.77350269}}
\put(0,0){\line(0,1){5.77350269}}
\put(0,5.77350269){\line(-1732,1000){5}}
\put(-10,5.77350269){\line(-1732,1000){5}}
\put(-10,5.77350269){\line(1732,1000){5}}
\put(0,5.77350269){\line(1732,1000){5}}
}
\multiput(10, -29.6413549)(10,0){3}{
\put(-10,0){\line(0,1){5.77350269}}
\put(0,0){\line(0,1){5.77350269}}
\put(0,5.77350269){\line(-1732,1000){5}}
\put(-10,5.77350269){\line(-1732,1000){5}}
\put(-10,5.77350269){\line(1732,1000){5}}
\put(0,5.77350269){\line(1732,1000){5}}
}
\multiput(15, -38.3016936)(10,0){2}{
\put(0,5.77350269){\line(-1732,1000){5}}
\put(-10,5.77350269){\line(-1732,1000){5}}
\put(-10,5.77350269){\line(1732,1000){5}}
\put(0,5.77350269){\line(1732,1000){5}}}

\put(15,-9.43392605){\makebox(0,0){$0$}}
\put(5,-9.43392605){\makebox(0,0){$5$}}
\put(25,-9.43392605){\makebox(0,0){$2$}}
\put(10,-0.77358733){\makebox(0,0){$6$}}
\put(20,-0.77358733){\makebox(0,0){$1$}}
\put(10,-18.0942648){\makebox(0,0){$4$}}
\put(20,-18.0942648){\makebox(0,0){$3$}}
%
\thinlines


\put(60,-16.43392605){\makebox(0,0){$\overrightarrow{\quad\xi\quad}$}}

\put(80,-38.3014395){\line(1,0){34.6410162}}
\put(80,-26.7544341){\line(1,0){46.1880215}}
\put(80,-15.2074287){\line(1,0){57.7350269}}
\put(80.5470054,-3.66042332){\line(1,0){57.7350269}}
\put(91.5470054,7.88658206){\line(1,0){46.1880215}}
\put(103.094011, 19.4335874){\line(1,0){34.6410162}}

\put(80,-38.3014395){\line(0,1){34.6410162}}
\put(91.5470054,-38.3014395){\line(0,1){46.1880215}}
\put(103.094011,-38.3014395){\line(0,1){57.7350269}}
\put(114.641016,-38.3014395){\line(0,1){57.7350269}}
\put(126.188022,-26.7544341){\line(0,1){46.1880215}}
\put(137.735027,-15.2074287){\line(0,1){34.6410162}}

\put(108.867513,-9.43392605){\makebox(0,0){$0$}}
\put( 97.3205081,-9.43392605){\makebox(0,0){$5$}}
\put(120.414519,-9.43392605){\makebox(0,0){$2$}}
\put(108.867513,2.11307933){\makebox(0,0){$6$}}
\put(97.3205081,-20.9809314){\makebox(0,0){$4$}}
\put(120.414519,2.11307933){\makebox(0,0){$1$}}
\put(108.867513,-20.9809314){\makebox(0,0){$3$}}
\end{picture}
\caption{The hexagonal model translation} \label{fig:hexmodel}
\end{figure}

We will often use an isomorphic representation of the hexagonal model
which will be of importance in the sequel. This representation
has $\Z^2$ as the set of vertices. Each point $(x,y)\in\Z^2$ has
the following neighboring vertices,
$$\{(x+a,y+b) ~|~ a,b\in\{-1,0,1\}, a+b\neq 0\}.$$ It may be shown
that the two models are isomorphic by using the mapping $\xi: \R^2
\rightarrow \R^2$, which is defined by
$\xi(x,y)=(x+\frac{y}{\sqrt{3}},\frac{2y}{\sqrt{3}})$. The effect of
the mapping on the neighbor set is shown in
Fig.~\ref{fig:hexmodel}. From now on, slightly changing notation, we
will also refer to this representation as the hexagonal model. Using
this new notation the neighbors of point $(i,j)$ are
\begin{multline*}
\{(i-1,j-1),(i-1,j),(i,j-1),(i,j+1),\hfill\\
\hfill (i+1,j),(i+1,j+1)\}.
\end{multline*}

The {\it hexagonal distance} $d(x,y)$ between two points $x$ and
$y$ in the hexagonal grid is the smallest $r$ such that there
exists a path with $r+1$ points $x=p_1,p_2,\ldots,p_{r+1}=y$,
where $p_i$ and $p_{i+1}$ are adjacent points in the hexagonal
grid.

\subsection{Distinct difference configurations}

We will now define our basic notation for the DDCs we will focus on.

\begin{definition}
\label{def:DDC}
A \emph{Euclidean square distinct difference configuration}
$\DD(m,r)$ is a set of $m$ dots placed in a square grid such that the
following two properties are satisfied:
\begin{enumerate}
\item Any two of the
dots in the configuration are at Euclidean distance at most $r$
apart.
\item All the $\binom{m}{2}$ differences between pairs of dots are
distinct either in length or in slope.
\end{enumerate}
\end{definition}

We will also study three more classes of DDCs: A \emph{square distinct
difference configuration} $\oDD(m,r)$ is defined by replacing
`Euclidean distance' by `Manhattan distance' in Definition~\ref{def:DDC};
a \emph{Euclidean hexagonal distinct difference configuration}
$\DD^*(m,r)$ is defined by replacing `square grid' by `hexagonal
grid' in Definition~\ref{def:DDC}; a \emph{hexagonal distinct
difference configuration} $\oDD^*(m,r)$ is defined by replacing `square
grid' by `hexagonal grid', and `Euclidean distance' by `hexagonal distance'
in Definition~\ref{def:DDC}.

%
%
%

In the application in~\cite{BEMP}, dots in the DDC are associated with
sensor nodes, and their position in the square or hexagonal grid
corresponds to a sensor's position. The parameter $r$ corresponds to a
sensor's wireless communication rage. So the most relevant distance
measure for the application we have in mind is the Euclidean
distance. Moreover, as the best packing of circles on a surface is by
arranging the circles in a hexagonal grid (see~\cite{CoSl93}), the
hexagonal model may be often be better from a practical point of
view. But the Manhattan and hexagonal distances are combinatorially
natural measures to consider, and our results for these distance
measures are sharper. Note that Manhattan and hexagonal distance are
both reasonable approximations to Euclidean distance (hexagonal
distance being the better approximation). Indeed, since the distinct
differences property does not depend on the distance measure used, it
is not difficult to show that a $\oDD(m,r)$ is a $\DD(m,r)$, and a
$\DD(m,r)$ is a $\oDD(m,\lceil\sqrt{2}r\rceil)$. Similarly a
$\oDD^*(m,r)$ is a $\DD^*(m,r)$, and a $\DD^*(m,r)$ is an
$\oDD^*(m,\lceil(2/\sqrt{3})r\rceil)$.

\subsection{Small parameters}

\begin{figure}
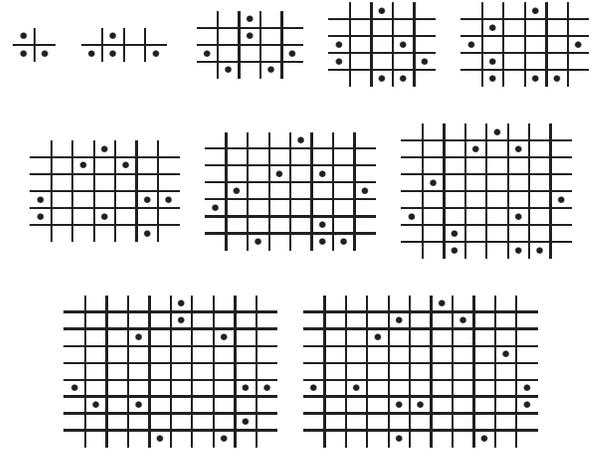

{\centering\tiny \setlength{\arraycolsep}{.4 \arraycolsep}
$\begin{array}{c|c}
\bullet & \\
\hline \bullet & \bullet
\end{array}$\quad\quad
$\begin{array}{c|c|c|c}
&\bullet&\phantom{\bullet}&\\
\hline \bullet&\bullet&&\bullet
\end{array}$ \quad\quad
$\begin{array}{c|c|c|c|c}
&&\bullet&&\\
\hline
&&\bullet&&\\
\hline
\bullet&&&&\bullet\\
\hline &\bullet&&\bullet&
\end{array}$\quad\quad
$\begin{array}{c|c|c|c|c}
&&\bullet&&\\
\hline
&\phantom{\bullet}&&&\\
\hline
\bullet&&&\bullet&\\
\hline
\bullet&&&&\bullet\\
\hline &&\bullet&\bullet&
\end{array}$\quad\quad
$\begin{array}{c|c|c|c|c|c}
&&\phantom{\bullet}&\bullet&&\\
\hline
&\bullet&&&&\\
\hline
\bullet&&&&&\bullet\\
\hline
&\bullet&&&&\\
\hline &\bullet&&\bullet&\bullet&
\end{array}$\\
\vspace{.5cm} $\begin{array}{c|c|c|c|c|c|c}
&&&\bullet&&&\\
\hline
&&\bullet&&\bullet&&\\
\hline
&\phantom{\bullet}&&&&&\\
\hline
\bullet&&&&&\bullet&\bullet\\
\hline
\bullet&&&\bullet&&&\\
\hline
&&&&&\bullet&\\
\end{array}$\quad\quad
$\begin{array}{c|c|c|c|c|c|c|c}
&&&&\bullet&&&\\
\hline
\phantom{\bullet}&&&&&&&\\
\hline
&&&\bullet&&\bullet&&\\
\hline
&\bullet&&&&&&\bullet\\
\hline
\bullet&&&&&&&\\
\hline
&&&&&\bullet&&\\
\hline &&\bullet&&&\bullet&\bullet&
\end{array}$\quad\quad
$\begin{array}{c|c|c|c|c|c|c|c}
&&&&\bullet&&&\\
\hline
&&&\bullet&&\bullet&&\\
\hline
\phantom{\bullet}&&&&&&&\\
\hline
&\bullet&&&&&&\\
\hline
&&&&&&&\bullet\\
\hline
\bullet&&&&&\bullet&&\\
\hline
&&\bullet&&&&&\\
\hline &&\bullet&&&\bullet&\bullet&
\end{array}$\\
\vspace{.5cm} $\begin{array}{c|c|c|c|c|c|c|c|c|c}
&&&&&\bullet&&&&\\
\hline
&&&&&\bullet&&&&\\
\hline
&&&\bullet&&&&\bullet&&\\
\hline
&&\phantom{\bullet}&&&&&&&\\
\hline
&&&&&&\phantom{\bullet}&&&\\
\hline
\bullet&&&&&&&&\bullet&\bullet\\
\hline
&\bullet&&\bullet&&&&&&\\
\hline
&&&&&&&&\bullet&\\
\hline &&&&\bullet&&&\bullet&&
\end{array}$\quad\quad
$\begin{array}{c|c|c|c|c|c|c|c|c|c|c}
&&&&&&\bullet&&&&\\
\hline
&&&&\bullet&&&\bullet&&&\\
\hline
&&&\bullet&&&&&&&\\
\hline
&&&&&&&&&\bullet&\\
\hline
&\phantom{\bullet}&&\phantom{\bullet}&&&&&&&\\
\hline
\bullet&&\bullet&&&&&&&&\bullet\\
\hline
&&&&\bullet&\bullet&&&&&\bullet\\
\hline
\phantom{\bullet}&&&&&&&&&&\\
\hline &&&&\bullet&&&&\bullet&&
\end{array}
$\\} \caption{Square distinct difference configurations with the largest
number of dots possible for $r=2,3,\dotsc,11$.}
\label{smallex}
\end{figure}

For small values of $r$, we used a backtrack search to exhaustively find a
$\oDD(m,r)$ with $m$ as large as possible. The search shows that for
$r=2,3$, the largest such $m$ are 3 and 4 respectively, and for
$4\leq r \leq 11$ the largest possible $m$ is
$r+2$. Fig.~\ref{smallex} contains examples of configurations meeting
those bounds.

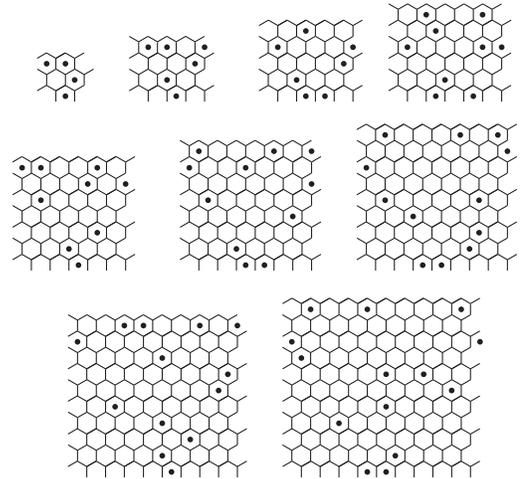
\begin{figure}[t]
\centering \setlength{\unitlength}{.25mm}
\begin{picture}(30,30)(-5,0)
\multiput(10,5)(10,0){2}{ \put(0,0){\line(0,1){5.77350269}}
\put(0,5.77350269){\line(-1732,1000){5}}
\put(-10,5.77350269){\line(1732,1000){5}}
\put(0,5.77350269){\line(1732,1000){5}} }
\multiput(5,13.66033872)(10,0){3}{ \put(0,0){\line(0,1){5.77350269}}
\put(0,5.77350269){\line(-1732,1000){5}}
\put(0,5.77350269){\line(1732,1000){5}} }
\multiput(10,22.3206774)(10,0){2}{ \put(0,0){\line(0,1){5.77350269}}
\put(0,5.77350269){\line(-1732,1000){5}}
\put(-10,5.77350269){\line(1732,1000){5}}
\put(0,5.77350269){\line(1732,1000){5}} }
\put(15,7.88675135){\circle*{3}} \put(20,16.5470901){\circle*{3}}
\put(15,25.2074288){\circle*{3}} \put(5,25.2074288){\circle*{3}}
\end{picture}
\quad
\begin{picture}(50,50)(-5,0)
\multiput(10,5)(10,0){4}{ \put(0,0){\line(0,1){5.77350269}}
\put(0,5.77350269){\line(-1732,1000){5}}
\put(-10,5.77350269){\line(1732,1000){5}}
\put(0,5.77350269){\line(1732,1000){5}} }
\multiput(5,13.66033872)(10,0){4}{ \put(0,0){\line(0,1){5.77350269}}
\put(0,5.77350269){\line(-1732,1000){5}}
\put(0,5.77350269){\line(1732,1000){5}} }
\multiput(10,22.3206774)(10,0){4}{ \put(0,0){\line(0,1){5.77350269}}
\put(0,5.77350269){\line(-1732,1000){5}}
\put(-10,5.77350269){\line(1732,1000){5}}
\put(0,5.77350269){\line(1732,1000){5}} }
\multiput(5,30.9810162)(10,0){4}{ \put(0,0){\line(0,1){5.77350269}}
\put(0,5.77350269){\line(-1732,1000){5}}
\put(0,5.77350269){\line(1732,1000){5}} }
\put(25,7.88675135){\circle*{3}} \put(20,16.5470901){\circle*{3}}
\put(35,25.2074288){\circle*{3}} \put(40,33.8677675){\circle*{3}}
\put(20,33.8677675){\circle*{3}} \put(10,33.8677675){\circle*{3}}

\end{picture}
\quad
\begin{picture}(50,50)(-5,0)
\multiput(0,0)(0,17.3206774){3} {\multiput(10,5)(10,0){5}{
\put(0,0){\line(0,1){5.77350269}}
\put(0,5.77350269){\line(-1732,1000){5}}
\put(-10,5.77350269){\line(1732,1000){5}}
\put(0,5.77350269){\line(1732,1000){5}} }}
\multiput(0,0)(0,17.3206774){2}{ \multiput(5,13.66033872)(10,0){5}{
\put(0,0){\line(0,1){5.77350269}}
\put(0,5.77350269){\line(-1732,1000){5}}
\put(0,5.77350269){\line(1732,1000){5}} }}
\put(25,7.88675135){\circle*{3}} \put(35,7.88675135){\circle*{3}}
\put(20,16.5470901){\circle*{3}} \put(45,25.2074288){\circle*{3}}
\put(50,33.8677675){\circle*{3}} \put(10,33.8677675){\circle*{3}}
\put(25,42.5281062){\circle*{3}}
\end{picture}
\quad
\begin{picture}(60,60)(-5,0)
\multiput(0,0)(0,17.3206774){3} {\multiput(10,5)(10,0){6}{
\put(0,0){\line(0,1){5.77350269}}
\put(0,5.77350269){\line(-1732,1000){5}}
\put(-10,5.77350269){\line(1732,1000){5}}
\put(0,5.77350269){\line(1732,1000){5}} }}
\multiput(0,0)(0,17.3206774){3}{ \multiput(5,13.66033872)(10,0){6}{
\put(0,0){\line(0,1){5.77350269}}
\put(0,5.77350269){\line(-1732,1000){5}}
\put(0,5.77350269){\line(1732,1000){5}} }}
\put(25,7.88675135){\circle*{3}} \put(45,7.88675135){\circle*{3}}
\put(30,16.5470901){\circle*{3}} \put(60,33.8677675){\circle*{3}}
\put(50,33.8677675){\circle*{3}} \put(10,33.8677675){\circle*{3}}
\put(25,42.5281062){\circle*{3}} \put(20,51.1884449){\circle*{3}}
\put(50,51.1884449){\circle*{3}}
\end{picture}\\

\begin{picture}(70,70)(-5,0)
\multiput(0,0)(0,17.3206774){4} {\multiput(10,5)(10,0){6}{
\put(0,0){\line(0,1){5.77350269}}
\put(0,5.77350269){\line(-1732,1000){5}}
\put(-10,5.77350269){\line(1732,1000){5}}
\put(0,5.77350269){\line(1732,1000){5}} }}
\multiput(0,0)(0,17.3206774){3}{ \multiput(5,13.66033872)(10,0){6}{
\put(0,0){\line(0,1){5.77350269}}
\put(0,5.77350269){\line(-1732,1000){5}}
\put(0,5.77350269){\line(1732,1000){5}} }}
\put(35,7.88675135){\circle*{3}} \put(30,16.5470901){\circle*{3}}
\put(45,25.2074288){\circle*{3}} \put(15,42.5281062){\circle*{3}}
\put(40,51.1884449){\circle*{3}} \put(60,51.1884449){\circle*{3}}
\put(45,59.8487837){\circle*{3}} \put(15,59.8487837){\circle*{3}}
\put(5,59.8487837){\circle*{3}}
\end{picture}
\quad
\begin{picture}(80,80)(-5,0)
\multiput(0,0)(0,17.3206774){4} {\multiput(10,5)(10,0){7}{
\put(0,0){\line(0,1){5.77350269}}
\put(0,5.77350269){\line(-1732,1000){5}}
\put(-10,5.77350269){\line(1732,1000){5}}
\put(0,5.77350269){\line(1732,1000){5}} }}
\multiput(0,0)(0,17.3206774){4}{ \multiput(5,13.66033872)(10,0){7}{
\put(0,0){\line(0,1){5.77350269}}
\put(0,5.77350269){\line(-1732,1000){5}}
\put(0,5.77350269){\line(1732,1000){5}} }}
\put(35,7.88675135){\circle*{3}} \put(45,7.88675135){\circle*{3}}
\put(30,16.5470901){\circle*{3}} \put(60,33.8677675){\circle*{3}}
\put(15,42.5281062){\circle*{3}} \put(70,51.1884449){\circle*{3}}
\put(35,59.8487837){\circle*{3}} \put(5,59.8487837){\circle*{3}}
\put(10,68.5091224){\circle*{3}} \put(50,68.5091224){\circle*{3}}
\put(70,68.5091224){\circle*{3}}
\end{picture}\quad
\begin{picture}(90,90)(-5,0)
\multiput(0,0)(0,17.3206774){5} {\multiput(10,5)(10,0){8}{
\put(0,0){\line(0,1){5.77350269}}
\put(0,5.77350269){\line(-1732,1000){5}}
\put(-10,5.77350269){\line(1732,1000){5}}
\put(0,5.77350269){\line(1732,1000){5}} }}
\multiput(0,0)(0,17.3206774){4}{ \multiput(5,13.66033872)(10,0){8}{
\put(0,0){\line(0,1){5.77350269}}
\put(0,5.77350269){\line(-1732,1000){5}}
\put(0,5.77350269){\line(1732,1000){5}} }}
\put(35,7.88675135){\circle*{3}} \put(45,7.88675135){\circle*{3}}
\put(60,16.5470901){\circle*{3}} \put(65,25.2074288){\circle*{3}}
\put(30,33.8677675){\circle*{3}} \put(15,42.5281062){\circle*{3}}
\put(65,42.5281062){\circle*{3}} \put(5,59.8487837){\circle*{3}}
\put(80,68.5091224){\circle*{3}} \put(15,77.1694611){\circle*{3}}
\put(55,77.1694611){\circle*{3}} \put(75,77.1694611){\circle*{3}}
\end{picture}\\
\begin{picture}(100,100)(-5,0)
\multiput(0,0)(0,17.3206774){5} {\multiput(10,5)(10,0){9}{
\put(0,0){\line(0,1){5.77350269}}
\put(0,5.77350269){\line(-1732,1000){5}}
\put(-10,5.77350269){\line(1732,1000){5}}
\put(0,5.77350269){\line(1732,1000){5}} }}
\multiput(0,0)(0,17.3206774){5}{ \multiput(5,13.66033872)(10,0){9}{
\put(0,0){\line(0,1){5.77350269}}
\put(0,5.77350269){\line(-1732,1000){5}}
\put(0,5.77350269){\line(1732,1000){5}} }}
\put(55,7.88675135){\circle*{3}} \put(50,16.5470901){\circle*{3}}
\put(65,25.2074288){\circle*{3}} \put(50,33.8677675){\circle*{3}}
\put(25,42.5281062){\circle*{3}} \put(80,51.1884449){\circle*{3}}
\put(85,59.8487837){\circle*{3}} \put(50,68.5091224){\circle*{3}}
\put(5,77.1694611){\circle*{3}} \put(30,85.8297998){\circle*{3}}
\put(40,85.8297998){\circle*{3}} \put(70,85.8297998){\circle*{3}}
\put(90,85.8297998){\circle*{3}}
\end{picture}\quad
\begin{picture}(100,110)(-5,0)
\multiput(0,0)(0,17.3206774){6} {\multiput(10,5)(10,0){10}{
\put(0,0){\line(0,1){5.77350269}}
\put(0,5.77350269){\line(-1732,1000){5}}
\put(-10,5.77350269){\line(1732,1000){5}}
\put(0,5.77350269){\line(1732,1000){5}} }}
\multiput(0,0)(0,17.3206774){5}{ \multiput(5,13.66033872)(10,0){10}{
\put(0,0){\line(0,1){5.77350269}}
\put(0,5.77350269){\line(-1732,1000){5}}
\put(0,5.77350269){\line(1732,1000){5}} }}
\put(45,7.88675135){\circle*{3}} \put(55,7.88675135){\circle*{3}}
\put(60,16.5470901){\circle*{3}} \put(30,33.8677675){\circle*{3}}
\put(55,42.5281062){\circle*{3}} \put(90,51.1884449){\circle*{3}}
\put(55,59.8487837){\circle*{3}} \put(75,59.8487837){\circle*{3}}
\put(10,68.5091224){\circle*{3}} \put(5,77.1694611){\circle*{3}}
\put(105,77.1694611){\circle*{3}} \put(95,94.4901385){\circle*{3}}
\put(15,94.4901385){\circle*{3}} \put(45,94.4901385){\circle*{3}}
\end{picture}
\caption{Hexagonal distinct difference configurations with the largest number
of dots possible for $r=2,3,\dotsc,10$.}
\label{fig:smallex_hex}
\end{figure}

Similarly, we found the best configurations $\oDD^*(m,r)$ in the hexagonal
grid (see Figure~\ref{fig:smallex_hex}) for $2\leq r\leq 10$.

\section{Anticodes and DDCs}
\label{sec:anticodes}

In this section we will show a trivial connection between DDCs and
maximal anticodes. This leads to a short investigation of maximal
anticodes in the square and the hexagonal grids. We find all maximal
anticodes in these two models under the Manhattan and hexagonal
distance measures respectively.

An {\it anticode} of diameter $r$ in the two-dimensional grid
(square or hexagonal) is a set $\cS$ of points such that for each
pair of points $x,y \in \cS$ we have $d(x,y) \leq r$, where the
distance can be Manhattan, hexagonal, or Euclidean. An anticode
$\cS$ of diameter $r$ is said to be {\it optimal} if there is no
anticode $\cS'$ of diameter $r$ such that $| \cS' | > | \cS |$. An
anticode $\cS$ of diameter $r$ is said to be {\it maximal} if $\{
x \} \cup \cS$ has diameter greater than $r$ for any $x \notin \cS$.
Anticodes are important structures in various aspects of coding
theory and extremal
combinatorics~\cite{AAK01,AhKh96,AhKh98,Del73,ESV06,MaZh95,ScEt02}.

The following two results provide an obvious connection between DDCs
and anticodes.
\begin{lemma}
\label{lem:extend_max} Any anticode $\cS$ of diameter $r$ is contained in
a maximal anticode $\cS'$ of diameter $r$.
\end{lemma}
\begin{corollary}
A $\DD(m,r)$ is contained in a maximal anticode of (Euclidean)
diameter~$r$. The same statement holds for a $\oDD(m,r)$, $\DD^*(m,r)$
or $\oDD^*(m,r)$ when the appropriate distance measure is used.
\end{corollary}

\subsection{Maximal anticodes in the square grid}

We start by defining three shapes in the square grid. We will
prove that these shapes are the only maximal anticodes in the square
grid when we use Manhattan distance.

\begin{figure}[tb]
\centering
\begin{tabular}{ccc}

\setlength{\unitlength}{.5mm}

\begin{picture}(25,35)(5,0)

\linethickness{.5 pt}
\put(15,5){\framebox(5,25){}}

\put(10,10){\framebox(15,15){}}

\put(5,15){\framebox(25,5){}}

\end{picture}&

\setlength{\unitlength}{.5mm}
\begin{picture}(30,40)(40,0)
\linethickness{.5 pt}
\put(40,15){\framebox(30,5){}}

\put(45,10){\framebox(20,15){}}

\put(50,5){\framebox(10,25){}}

\put(55,5){\line(0,1){25}}

\end{picture}&

\setlength{\unitlength}{.5mm}
\begin{picture}(30,40)(80,0)
\linethickness{.5 pt}
\put(80,15){\framebox(30,10){}}

\put(85,10){\framebox(20,20){}}

\put(90,5){\framebox(10,30){}}

\put(95,5){\line(0,1){30}}

\put(80,20){\line(1,0){30}}

\end{picture}\\
(a) & (b) & (c)
\end{tabular}

\caption{Maximal anticodes in the square grid}\label{genleespheres}

\end{figure}
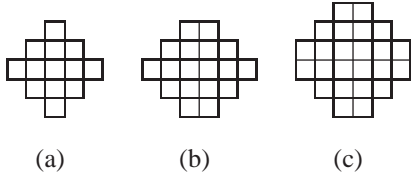

A {\it Lee sphere} of radius R is the shape in the square model
which consists of one point as centre and all positions of Manhattan
distance at most $R$ from this centre. The area of this Lee sphere
is $2R^2+2R+1$.  For the seminal paper on Lee spheres
see~\cite{GoWe70}. Figure~\ref{genleespheres}a illustrates a Lee sphere
of radius 2.

A {\it bicentred Lee sphere} of radius R is the shape in the square
model which consists of two centre points (a $2 \times 1$ or an $1
\times 2$ rectangle) and all positions of Manhattan distance at most
$R$ from at least one point of this centre. The area of this bicentred
Lee sphere is $2R^2+4R+2$. These shapes were used for
two-dimensional burst-correction in~\cite{BBV98}.
Figure~\ref{genleespheres}b illustrates a bicentred Lee sphere of
radius 2.

A {\it quadricentred Lee sphere} of radius R is the shape in the
square model which consists of four centre points (a $2 \times 2$
square) and all positions of Manhattan distance at most $R-1$ from at
least one point of this centre. The area of this quadricentred Lee
sphere is $2R^2+2R$. These shapes were defined using the name
`generalized Lee sphere' in~\cite{Etz03}.  Figure~\ref{genleespheres}c
illustrates a quadricentred Lee sphere of radius 3.

\begin{theorem}
\label{thm:max_sqr} $~$

\begin{itemize}
\item
For even $r$ there are two different types of maximal anticodes of
diameter $r$ in square grid: the Lee
sphere of radius $\frac{r}{2}$ and the quadricentred Lee sphere of
radius $\frac{r}{2}$.

\item For odd $r$ there is exactly one type of maximal anticode of
diameter $r$ in the square grid: the bicentred
Lee sphere of radius $\frac{r-1}{2}$.
\end{itemize}
\end{theorem}
\begin{IEEEproof}
Let $\cA$ be a maximal anticode of diameter $r$ in the square
grid.

Assume first that $r$ is even, so $r=2\rho$. We will embed $\cA$ in
the two-dimensional square grid in such a way that there is position
in $\cA$ on the line $y=x$, but no position below it. The Manhattan
distance between a point on the line $y=x+2\rho$ and a point on the
line $y=x$ is at least $2\rho$ and hence $\cA$ is bounded by the lines
$y=x$ and $y=x+2\rho$. Similarly, without loss of generality we can
assume that there is a position in $\cA$ on the line $y=-x$ or
$y=-x+1$ and no position below this line, so $\cA$ is bounded
by the lines $y=-x$ and $y=-x+2\rho$, or by the lines
$y=-x+1$ and $y=-x+2\rho+1$. These four lines define a Lee sphere of
radius $\rho$ or a quadricentred Lee sphere of radius $\rho$.

Now, assume that $r$ is odd, so $r=2\rho+1$. We will embed $\cA$ in
the two-dimensional square grid in a way that a position of $\cA$ lies
on the line $y=x$, but no position lies below it. The Manhattan
distance between a point on the line $y=x+2\rho+1$ and a point on the
line $y=x$ is at least $2\rho+1$ and hence $\cA$ is bounded by the
lines $y=x$ and $y=x+2\rho+1$. Similarly, without loss of generality
we can assume that there is a position of $\cA$ on the line $y=-x$ or
$y=-x+1$ and no position below this line, and so $\cA$ is bounded by
the lines $y=-x$ and $y=-x+2\rho+1$ or by the lines $y=-x+1$ and
$y=-x+2\rho+2$. In either case, these four lines define a bicentred
Lee sphere of radius $\rho$.
\end{IEEEproof}

Finally, the following theorem is interesting from a theoretical
point of view.
\begin{theorem}
There exists a $\oDD(m,r)$ for which the only maximal anticode of diameter $r$
containing it is a Lee sphere (bicentred Lee sphere, quadricentred Lee
sphere) of diameter $r$.
\end{theorem}
\begin{IEEEproof}
We provide the configurations that are needed. All the claims in the
proof below are readily verified and left to the reader.

When $r$ is odd, we may take two points on the same horizontal line
such that $d(x,y)=r$: this pair of points is in a bicentred Lee sphere
of radius $\frac{r-1}{2}$. When $r$ is even, the same example is
contained in a Lee sphere of radius $r/2$, but is not contained in a
quadricentred Lee sphere of radius $r/2$.

Let $r$ be even, and set $R=r/2$. The points $(0,R-1)$, $(0,R)$,
$(2R-2,0)$, $(2R-2,2R-1)$, and $(2R-1,R)$ form $\oDD(5,2R)$. This set
of points is not contained in a Lee sphere of radius $R$, but is
contained in a quadricentred Lee sphere of radius $R$.
\end{IEEEproof}

\subsection{Maximal anticodes in the hexagonal grid}

\begin{theorem}
There are exactly $\lceil \frac{r+1}{2} \rceil$ different types of
maximal anticodes of diameter $r$ in the hexagonal grid, namely the
anticodes $\cA_0,\cA_1,\ldots ,\cA_{\lceil \frac{r-1}{2}\rceil}$
defined in the proof below.
\end{theorem}
\begin{IEEEproof}
We consider the translation of the hexagonal grid into the square
grid. By shifting it appropriately, any maximal anticode $\cA$ of
diameter $r$ can be located inside an $(r+1) \times (r+1)$ square
$\cB$ with corners at $(0,0)$, $(0,r)$, $(r,0)$, and $(r,r)$.  Let $i$
be defined by the property that the lines $y=x-i$ contains a point of
$\cA$, but no point of $\cA$ lies below this line.

We claim that $i \geq 0$. To see this, assume for a contradiction that
$i < 0$. Then $\cA$ is contained in the region of $\cB$ bounded by
the lines $y=x-i$, $y=r$, and $x=0$. But the point $(0,r+1)$ outside
$\cB$ is within distance $r$ from all the points of this region,
contradicting the fact that $\cA$ is a maximal anticode. Thus, $i \geq
0$ and our claim follows.

All the points on the line $y=x-i$ that are inside $\cB$ are within
hexagonal distance $r$ from all points on the lines $y=x-i+j$, $0 \leq
j \leq r$, that lie inside $\cB$. All the points on the line $y=x-i$
inside $\cB$ have hexagonal distance greater than $r$ from all the
points on the line $y=x-i+r+1$. Hence, as $\cA$ is maximal, $\cA$
consists of all the points bounded by the lines $y=x-i$ and $y=x-i+r$
inside $\cB$. It is easy to verify that each one of the $r+1$
anticodes $\cA_i$ defined in this way is a maximal anticode. One can
readily verify that $\cA_i$ and $\cA_{r-i}$ are equivalent anticodes,
since $\cA_{r-i}$ is obtained by rotating $\cA_i$ by 180 degrees. So
the theorem follows.
\end{IEEEproof}

\begin{theorem}
Let $i$ be fixed, where $0 \leq i \leq \lceil \frac{r-1}{2} \rceil$.
There exists a $\oDD^*(m,r)$ for which the only maximal anticode of diameter $r$
containing it is of the form $\cA_i$.
\end{theorem}
\begin{IEEEproof}
Again, we provide the configurations, and leave the verification of
the details to the reader.

For each $i$, $1 \leq i \leq \frac{r-1}{2}$, $\cA_i$ has six corner
points. If we assign a dot to each corner point then we will obtain a
DDC which cannot be inscribed in another maximal anticode. When $r=2R$
these six points do not define a DDC in $\cA_R$. In this case we
assign seven dots to $\cA_R$ as follows. In four consecutive corner
points we assign a dot; in the next corner point we assign two dots in
the adjacent points on the boundary of $\cA_R$; in the last corner
point we assign a dot in the adjacent point on the boundary of $\cA_R$
towards the first corner point. When $i=0$, $\cA_i$ is a triangle and
has three corner points. If we assign a dot to each of these corner
points we will obtain a DDC which cannot be inscribed in another
maximal anticode.
\end{IEEEproof}

\begin{figure}[t]
\centering
\setlength{\unitlength}{.4mm}
\begin{picture}(50,60)(-10,-35)
\multiput(10,5)(10,0){3}{
\put(-10,0){\line(0,1){5.77350269}}
\put(0,0){\line(0,1){5.77350269}}
\put(0,5.77350269){\line(-1732,1000){5}}
\put(-10,5.77350269){\line(1732,1000){5}}
}
\multiput(5,-3.66033872)(10,0){4}{
\put(-10,0){\line(0,1){5.77350269}}
\put(0,0){\line(0,1){5.77350269}}
\put(0,5.77350269){\line(-1732,1000){5}}
\put(-10,5.77350269){\line(1732,1000){5}}
}
\multiput(0,-12.3206774)(10,0){5}{
\put(-10,0){\line(0,1){5.77350269}}
\put(0,0){\line(0,1){5.77350269}}
\put(0,5.77350269){\line(-1732,1000){5}}
\put(-10,5.77350269){\line(1732,1000){5}}
}
\multiput(5, -20.9810162)(10,0){4}{
\put(-10,0){\line(0,1){5.77350269}}
\put(0,0){\line(0,1){5.77350269}}
\put(0,5.77350269){\line(-1732,1000){5}}
\put(-10,5.77350269){\line(-1732,1000){5}}
\put(-10,5.77350269){\line(1732,1000){5}}
\put(0,5.77350269){\line(1732,1000){5}}
}
\multiput(10, -29.6413549)(10,0){3}{
\put(-10,0){\line(0,1){5.77350269}}
\put(0,0){\line(0,1){5.77350269}}
\put(0,5.77350269){\line(-1732,1000){5}}
\put(-10,5.77350269){\line(-1732,1000){5}}
\put(-10,5.77350269){\line(1732,1000){5}}
\put(0,5.77350269){\line(1732,1000){5}}
}
\multiput(15, -38.3016936)(10,0){2}{
\put(0,5.77350269){\line(-1732,1000){5}}
\put(-10,5.77350269){\line(-1732,1000){5}}
\put(-10,5.77350269){\line(1732,1000){5}}
\put(0,5.77350269){\line(1732,1000){5}}}
\end{picture}
\caption{Hexagonal sphere of radius 2}
\label{hexsphere}
\end{figure}
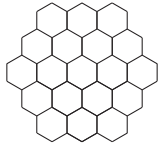

We now consider some basic properties of these $\lceil \frac{r+1}{2}
\rceil$ anticodes. First, the number of grid points in $\cA_i$ is
$(r+1)^2 - \frac{i(i+1)}{2} - \frac{(r-i)(r+1-i)}{2}=
\frac{(r+1)(r+2)}{2} +i(r-i)$. The smallest anticode is $\cA_0$, an
isosceles right triangle with base and height of length $r+1$
containing $\frac{(r+1)(r+2)}{2}$ points. The largest anticode is the
\emph{hexagonal sphere $\cA_{\lceil \frac{r-1}{2} \rceil}$ of radius
$r/2$}. The hexagonal sphere contains $\frac{3(r+1)^2}{4}$ points when
$r$ is odd, and contains $\frac{3r^2+6r+4}{4} = 3 \left( \frac{r}{2}
\right)^2 +3\left(\frac{r}{2}\right) +1$ points when $r$ is even. The hexagonal
sphere of radius R is the shape in the hexagonal model which
consists of a centre point and all positions in hexagonal distance at
most $R$ from this centre (Fig.~\ref{hexsphere}).

\subsection{Maximal anticodes with Euclidean distance}

It seems much more difficult to classify the maximal anticodes in the
square and hexagonal grids when we use Euclidean distance.  Note that
the representation of the hexagonal grid in the square grid does not
preserve Euclidean distances, and so we cannot use the map $\xi$. We
expect that the overall shape of a maximal anticode in both models
should be similar, since a maximal anticode in both models is just the
intersection of a maximal anticode in $\R^2$ with the centres of our
squares or hexagons respectively. But the `local' structure of an
anticode will be different: for example, in the hexagonal grid we can
have three dots that are pairwise at distance $r$, but this is not
possible in the square grid.

Because maximal anticodes in $\R^2$ determine the shape of maximal
anticodes in the square or hexagonal models, we conclude this section with
a brief description of such anticodes.

\begin{figure}\centering
\begin{tabular}{ccc}
\begin{picture}(50,70)(-25,-35)
\put(-20,0){\arc[-60,60]{40}}
\put(20,0){\arc[120,240]{40}}
\put(-20,0){\circle*{2.5}}
\put(20,0){\circle*{2.5}}
\linethickness{0.2 pt}
\put(-20,0){\line(1,0){40}}
\put(-2,1){\makebox(4,5){$r$}}
\end{picture}
&
\begin{picture}(50,70)(-25,-35)
\put(-20,0){\arc[-60,60]{40}}
\put(20,0){\arc[120,240]{40}}
\put(-20,0){\circle*{2.5}}
\put(20,0){\circle*{2.5}}
\linethickness{0.2 pt}
\put(-20,0){\line(1,0){40}}
\put(-2,1){\makebox(4,5){$r$}}
\linethickness{1 pt}
\put(0,0){\circle{40}}
\end{picture}
&
\begin{picture}(50,70)(-25,-35)
\put(-20,0){\arc[-60,60]{40}}
\put(20,0){\arc[120,240]{40}}
\put(-20,0){\circle*{2.5}}
\put(20,0){\circle*{2.5}}
\linethickness{0.2 pt}
\put(-20,0){\line(1,0){40}}
\put(-20,0){\line(1000,1732){20}}
\put(20,0){\line(-1000,1732){20}}
\put(-2,1){\makebox(4,5){$r$}}
\linethickness{1 pt}
\put(-20,0){\arc[0,60]{40}}
\put(20,0){\arc[120,180]{40}}
\put(-20,0){\arc[0,60]{40}}
\put(0,34.64){\arc[240,300]{40}}
\put(0,34.63){\circle*{2.5}}
\end{picture}\\
(a) & (b) & (c)
\end{tabular}
\caption{Anticodes in the Euclidean distance}\label{fig:anti}
\end{figure}
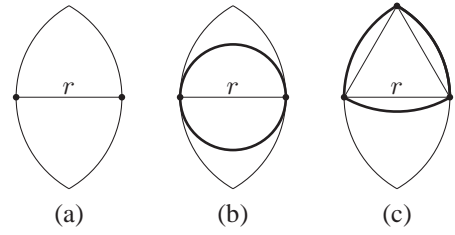

An anticode is confined to the area as
depicted in Figure~\ref{fig:anti}a, where dots are two
elements in the anticode at distance $r$. The most obvious maximal
anticode is a circle of diameter $r$ depicted in
Figure~\ref{fig:anti}b. Another maximal anticode is
depicted in Figure~\ref{fig:anti}c, and is constructed by taking
three dots at the vertices of an equilateral triangle of side $r$, and
intersecting the circles of radius $r$ about these dots.  Between the
`triangular' anticode and the circle there are infinitely many other
maximal anticodes. We will need the following `isoperimetrical'
theorem; see Littlewood~\cite[Page~32]{Lit86} for a proof.

\begin{theorem}
\label{thm:plane_anticode}
Let $\cA$ be a region of $\R^2$ of diameter $r$ and area $a$. Then
$a\leq (\pi/4) r^2$.
\end{theorem}

We remark that the example of a circle of diameter $r$ shows that the
bound of this theorem is tight.

\section{Upper bounds on the Number of Dots}
\label{sec:upper_bounds}

In this section we will provide asymptotic upper bounds on the number
of dots that can be contained in a DDC, using a technique due to Erd\H
os and Tur\'an~\cite{EGRT92,Erd41}. We start by considering upper
bounds on the number $m$ of dots in a $\oDD(m,r)$ and a $\oDD^*(m,r)$,
and then consider upper bounds in a $\DD(m,r)$ and $\DD^*(m,r)$. The
results for small parameters in Section~\ref{sec:models} might suggest
that a $\oDD(m,r)$ can always contain $r+2$ dots: our result
(Theorem~\ref{thm:up_Lee}) that $m\leq \frac{1}{\sqrt{2}}r+o(r)$
surprised us. Our techniques easily generalise to DDCs where we
restrict the dots to lie in various shapes in the grid not necessarily
related to distance measures: we end the section with a brief
discussion of this general situation.

\subsection{Manhattan and hexagonal distances}
\label{sec:upper}

\begin{lemma}
\label{anticode_coverage_lemma} Let $r$ be a non-negative integer. Let ${\cal A}$ be an anticode of
Manhattan diameter $r$ in the square grid. Let $\ell$ be a positive
integer such that $\ell\leq r$, and let $w$ be the number of Lee
spheres of radius $\ell$ that intersect ${\cal A}$ non-trivially. Then
$w\leq\frac{1}{2}(r+2\ell)^2+O(r)$.
\end{lemma}
\begin{IEEEproof}
Let $\cal A'$ be the set of centres of the Lee spheres we are
considering, so $w=|\cA'|$. We claim that $\cA'$ is an anticode of
diameter at most $r+2\ell$. To see this, let $c,c'\in\cA'$. Since the
spheres of radius $\ell$ about $c$ and $c'$ intersect $\cA$
non-trivially, there exist elements $a,a'\in\cA$ such that
$d(c,a)\leq\ell$ and $d(c',a')\leq\ell$. But then
\[
d(c,c')\leq d(c,a)+d(a,a')+d(a',c')\leq \ell+r+\ell=r+2\ell,
\]
and so our claim follows.

Let $\cA''$ be a maximal anticode of diameter $r+2\ell$ containing
$\cA'$. Theorem~\ref{thm:max_sqr} implies that $\cA''$ is a Lee
sphere, bicentred Lee sphere or quadricentred Lee sphere of radius
$R$, where $R=\lfloor (r+2\ell)/2\lfloor$. In all three cases,
$|\cA''|=2R^2+O(R)=\frac{1}{2}(r+2\ell)^2+O(r)$. Since
\[
w=|\cA'|\leq|\cA''|,
\]
the lemma follows.
\end{IEEEproof}

\begin{theorem}
\label{thm:up_Lee} If a $\oDD(m,r)$ exists, then
\[
m\leq
\tfrac{1}{\sqrt{2}} r + (3/2^{4/3})r^{2/3}+O(r^{1/3}).
\]
\end{theorem}
\begin{IEEEproof}
We begin by giving a simple argument that leads to a linear bound
on $m$ in terms of $r$, with an inferior leading term to the bound
in the statement of the theorem. There are $2r^2+2r$ non-zero
vectors of Manhattan length $r$ or less, where a vector is
a line with direction which connects two points. The distinct
difference property implies that each such vector arises at most
once as the vector difference of a pair of dots. Since a
configuration of $m$ dots gives rise to $m(m-1)$ vector
differences, we find that
\[
m(m-1)\leq 2r^2+2r.
\]
In particular, we see that $m\leq \sqrt{2}r+o(r)=O(r)$.

We now establish the bound of the theorem. Since all the dots are at
distance at most $r$, we see that all dots are contained in a fixed
anticode $\cal A$ of the square grid of diameter $r$. Set
$\ell=\lfloor \alpha r^{2/3}\rfloor$, where we will choose the
constant $\alpha$ later so as to optimize our bound.  We cover
$\cal A$ with all the `small' Lee spheres of radius $\ell$ that
intersect $\cal A$ nontrivially. Every point of $\cal A$ is contained
in exactly $a$ small Lee spheres, where $a=2\ell^2+2\ell+1$.  Moreover, by
Lemma~\ref{anticode_coverage_lemma}, we have used $w$ small Lee
spheres, where $w\leq \frac{1}{2}(r+2\ell)^2+O(r)$.

Let $m_i$ be the number of dots in the $i$th small Lee sphere. Let
$\mu$ be the mean of the integers $m_i$. Since every dot is
contained in exactly $a$ small Lee spheres,
$\mu=am/w$. We aim to show that
\begin{equation}
\label{eqn:key_inequality} w(\mu^2-\mu)\leq \sum_{i=1}^w
m_i(m_i-1) \leq a(a-1).
\end{equation}

The first inequality in~\eqref{eqn:key_inequality} follows from
expanding the non-negative sum $\sum_{i=1}^w (\mu-m_i)^2$, so it
remains to show the second inequality.

The sum $\sum_{i=1}^w m_i(m_i-1)$ counts the number of pairs
$({\cal L},d)$ where ${\cal L}$ is a small Lee sphere and $d$ is a
vector difference between two dots in $\ell$. Every difference $d$
arises from a unique ordered pair of dots in $\cal A$, since the dots form
a distinct difference configuration. Thus
\[
\sum_{i=1}^w m_i(m_i-1)\leq \sum_d k(d),
\]
where we sum over all non-zero vector differences $d$ and where
$k(d)$ is the number of Lee spheres of radius $\ell$ that contain
any fixed pair of dots with vector difference $d$. If we assume
that the first element of the pair of dots with vector difference
$d$ lies at the origin, we see that
\[
\sum_d k(d) = a(a-1)
\]
since there are exactly $a$ Lee spheres of radius
$\ell$ containing the origin, and each such sphere contributes $1$
to $k(d)$ for exactly $a-1$ values of $d$. Thus we have
established~\eqref{eqn:key_inequality}.

Now, the inequality~\eqref{eqn:key_inequality} together with the
fact that $\mu=am/w$ imply that
\[
(\mu-1)m\leq a-1\leq a,
\]
and so
\begin{equation}
\label{square_eqn}
m^2\leq w\left(1+\frac{m}{a}\right).
\end{equation}
By Lemma~\ref{anticode_coverage_lemma},
\[
\sqrt{w}\leq \frac{1}{\sqrt{2}}r\left(1+\frac{2\ell}{r}+O(r^{-1})\right),
\]
and we have that
\[
\sqrt{(1+(m/a))}=1+m/(2a)+O((m/a)^2).
\]
Since $m=O(r)$ and $a\geq 2\alpha^2r^{4/3}$, these two inequalities
combine with~\eqref{square_eqn} to show that
\begin{equation}
\label{square_upp_eqn}
m\leq \frac{1}{\sqrt{2}}r\left(1+2\alpha
r^{-1/3}+\frac{m}{4\alpha^2r^{4/3}}+O(r^{-2/3})\right).
\end{equation}
Since $m=O(r)$, this inequality implies that $m\leq
\frac{1}{\sqrt{2}}r+O(r^{2/3})$. Combining this tighter bound
with~\eqref{square_upp_eqn} we find that
\[
m\leq
\frac{1}{\sqrt{2}}r\left(1+\left(2\alpha+\frac{1}{4\sqrt{2}\alpha^2}\right)r^{-1/3}+O(r^{-2/3})\right).
\]
The expression $2\alpha+1/(4\sqrt{2}\alpha^2)$ is minimized when
$\alpha=2^{-5/6}$ at the value $3/2^{5/6}$, so choosing this value for
$\alpha$ we deduce that
\begin{eqnarray*}
m&\leq\frac{1}{\sqrt{2}}r\left(1+\frac{3}{2^{5/6}}r^{-1/3}+O(r^{-2/3})\right)\\
&=\frac{1}{\sqrt{2}}r+\frac{3}{2^{4/3}}r^{2/3}+O(r^{1/3}),
\end{eqnarray*}
as required.
\end{IEEEproof}

We now look at the hexagonal grid.

\begin{lemma}
\label{hex_anticode_coverage_lemma} Let $r$ be a non-negative
integer. Let ${\cal A}$ be an anticode of hexagonal diameter $r$ in
the hexagonal grid. Let $\ell$ be a positive integer such that
$\ell\leq r$, and let $w$ be the number of hexagonal spheres of radius
$\ell$ that intersect ${\cal A}$ non-trivially. Then $w\leq
\frac{3}{4}(r+2\ell)^2+O(r)$.
\end{lemma}
\begin{IEEEproof}
The set of centres of the hexagonal spheres of radius $\ell$ that
have non-trivial intersection with $\cal A$ clearly form an
anticode of diameter at most $r+2\ell$. Therefore the number $w$
of such spheres is bounded by the maximal size of such an
anticode. The results on the maximal anticodes in the hexagonal
metric in Section~\ref{sec:anticodes} imply that
\begin{align*}
w&\leq\tfrac{1}{4}(3(r+2\ell)^2+6(r+2\ell)+4)\\
&=\tfrac{3}{4}(r+2\ell)^2+O(r),
\end{align*}
as required.
\end{IEEEproof}

\begin{theorem}
\label{thm:up_hex} If a $\oDD^*(m,r)$ exists, then
\[
m\leq
\tfrac{\sqrt{3}}{2} r + (3^{4/3}2^{-5/3})r^{2/3}+O(r^{1/3}).
\]
\end{theorem}
\begin{IEEEproof}
The dots in a $\oDD^*(m,r)$ form an anticode of diameter $r$. Let
$\ell=\lfloor 2^{-2/3}3^{-1/6}r^{2/3}\rfloor$. We may cover these
dots with the $w$ hexagonal spheres of radius $\ell$ that contain
one or more of these dots. By
Lemma~\ref{hex_anticode_coverage_lemma}, we have that $w\leq
\frac{3}{4}(r+2\ell)^2+O(r)$.

Using the fact that a hexagonal sphere of radius $\ell$ contains a
points in the hexagonal grid, where $a=3\ell^2+3\ell+1$, we may argue
exactly as in Theorem~\ref{thm:up_Lee} to produce the
bound~\eqref{square_eqn}. There are $O(r^2)$ vectors in the hexagonal
grid of hexagonal length $r$ or less, so the argument in
the first paragraph of Theorem~\ref{thm:up_Lee} shows that $m=O(r)$. Since $m/a=O(r^{-1/3})$, the bound~\eqref{square_eqn} implies that
\[
m\leq \sqrt{w}+O(r^{2/3}) = \frac{\sqrt{3}}{2}r+O(r^{2/3}).
\]
This bound on $m$ implies that
$m/a=2^{1/3}3^{-1/6}r^{-1/3}+O(r^{-2/3})$, and so
applying~\eqref{square_eqn} once more we obtain the bound of the
theorem, as required.
\end{IEEEproof}

One consequence of Theorem~\ref{thm:up_hex} is an answer to the ninth
question asked by Golomb and Taylor~\cite{GoTa84}: a {\it honeycomb
array} is a DDC in the hexagonal grid whose dots, when represented in
the square grid, form an $m\times m$ Costas array whose dots lie in
$m$ consecutive `North-East' diagonals. Honeycomb arrays are the
natural hexagonal analogue of Costas arrays. Do honeycomb arrays exist
for infinitely many $m$? The conjecture in~\cite{GoTa84} is that the
answer is YES. However, the answer is in fact NO, as the following
corollary to Theorem~\ref{thm:up_hex} shows.

\begin{corollary}
Honeycomb arrays exist for only a finite number of values of $m$.
\end{corollary}
\begin{IEEEproof}
The dots in a honeycomb array are contained in an anticode of diameter
at most $m-1$ (using hexagonal distance). Hence a honeycomb array is a
$\oDD^*(m,m-1)$.  But Theorem~\ref{thm:up_hex} shows that $m\leq
\frac{\sqrt{3}}{2}m+O(m^{2/3})$. Since $\frac{\sqrt{3}}{2}<1$, no
honeycomb array exists when $m$ is sufficiently large.
\end{IEEEproof}
In fact, numerical computations indicate that no honeycomb arrays
exist for $m\geq 1289$: for $m$ in this range, there is a suitable
choice of $\ell$ such that a honeycomb array violates the
inequality~\eqref{square_eqn}.

\subsection{Euclidean distance}

We now turn our attention to Euclidean distance. Our first lemma is
closely related to Gauss's circle problem:

\begin{lemma}
\label{lem:Euc_circ_square}
Let $\ell$ be a positive integer, and let $\cS$ be a (Euclidean)
circle of radius $\ell$ in the plane. Then the number of points of the
square grid contained in $\cS$ is $\pi \ell^2+O(\ell)$.
\end{lemma}
\begin{IEEEproof}
Let $c$ be the centre of $\cS$. Let $X$ be the set of points of the
square grid contained in $\cS$. Define $\cX$ to be the union of all
unit squares whose centres lie in $X$. Clearly $\cX$ has area $|X|$. The
maximum distance from the centre of a unit square to any point in the
unit square is at most $1/\sqrt{2}$, and so $\cX$ is contained in the
circle of radius $\ell+(1/\sqrt{2})$ with centre $c$. Similarly, every
point in a circle of radius $\ell-(1/\sqrt{2})$ with centre $c$ is
contained in $\cX$. Hence
\[
\pi(\ell-(1/\sqrt{2}))^2\leq |X|\leq\pi(\ell+(1/\sqrt{2}))^2,
\]
and so the lemma follows.
\end{IEEEproof}

\begin{lemma}
\label{lem:Euc_anti_square}
Let $r$ be a non-negative integer. Let $\cA$ be an anticode in the
square grid of Euclidean diameter $r$. Let $\ell$ be a positive
integer such that $\ell\leq r$, and let $w$ be the number of circles
of radius $\ell$ whose centres lie in the square grid and that intersect
$\cA$ non-trivially. Then $w\leq (\pi/4)(r+2\ell)^2+O(r)$.
\end{lemma}
\begin{IEEEproof}
As in Lemma~\ref{anticode_coverage_lemma}, it is not difficult to
see that the set $\cA'$ of centres of circles we are
considering form an anticode in the square grid of diameter at most
$r+2\ell$. Note that $w=|\cA'|$. Let $\cX$ be the union of the unit
squares whose centres lie in $\cA'$, so $\cX$ has area~$w$. The
maximum distance between the centre of a unit square and any other
point in this square is $1/\sqrt{2}$, and so $\cX$ is an anticode in
$\R^2$ of diameter at most $r+2\ell+(1/\sqrt{2})$. Hence, by
Theorem~\ref{thm:plane_anticode}, $w\leq
(\pi/4)(r+2\ell+(1/\sqrt{2}))^2$ and the lemma follows.
\end{IEEEproof}

\begin{theorem}
\label{thm:Euc_upper_square}
If a $\DD(m,r)$ exists, then
\[
m\leq
\tfrac{\sqrt{\pi}}{2}r+\tfrac{3\,\pi^{1/3}}{2^{5/3}}r^{2/3}+O(r^{1/3}).
\]
\end{theorem}
\begin{IEEEproof}
The proof is essentially the same as the proof of
Theorem~\ref{thm:up_hex}, using Lemma~\ref{lem:Euc_circ_square} to
bound the number $a$ of points in a sphere of radius $\ell$, and
using Lemma~\ref{lem:Euc_anti_square} instead of
Lemma~\ref{hex_anticode_coverage_lemma}. The bound of the theorem is
obtained if we set $\ell=\lfloor 1/(2^{2/3}\pi^{1/6})r^{2/3}\rfloor$.
\end{IEEEproof}

\begin{lemma}
\label{lem:Euc_circ_hex}
Let $\ell$ be a positive integer, and let $\cS$ be a (Euclidean)
circle of radius $\ell$ in the plane. Then the number of points of the
hexagonal grid contained in $\cS$ is $(2\pi/\sqrt{3}) \ell^2+O(\ell)$.
\end{lemma}
\begin{IEEEproof}
The proof of the lemma is essentially the same as the proof of
Lemma~\ref{lem:Euc_circ_square}. The hexagons whose centres form the
hexagonal grid have area $\sqrt{3}/2$, and the maximum distance of the
centre of a hexagon to any point in the hexagon is
$1/\sqrt{3}$. Define $X$ to be the set of points of the hexagonal grid
contained in $\cS$, and let $\cX$ be the union of all hexagons in our
grid whose centres lie in $X$. Clearly $\cX$ has area
$(\sqrt{3}/2)|X|$. The argument of Lemma~\ref{lem:Euc_circ_square} shows that
\[
\pi(\ell-(1/\sqrt{3}))^2\leq (\sqrt{3}/2)|X|\leq\pi(\ell+(1/\sqrt{3}))^2,
\]
and so the lemma follows.
\end{IEEEproof}

\begin{lemma}
\label{lem:Euc_anti_hex}
Let $r$ be a non-negative integer. Let $\cA$ be an anticode in the
square grid of Euclidean diameter $r$. Let $\ell$ be a positive
integer such that $\ell\leq r$, and let $w$ be the number of circles
of radius $\ell$ whose centres lie in the hexagonal grid and that intersect
$\cA$ non-trivially. Then $w\leq (\pi/(2\sqrt{3}))(r+2\ell)^2+O(r)$.
\end{lemma}
\begin{IEEEproof}
The proof of this lemma is essentially the same as the proof of
Lemma~\ref{lem:Euc_anti_square}. The argument there with appropriate
modifications shows that $(\sqrt{3}/2)w\leq
(\pi/4)(r+2\ell+(1/\sqrt{3}))^2$ (where the factor of $\sqrt{3}/2$
comes from the fact that the hexagons associated with our grid have area
$\sqrt{3}/2$).
\end{IEEEproof}

\begin{theorem}
\label{thm:Euc_upper_hex}
If a $\DD^*(m,r)$ exists, then
\[
m\leq
\tfrac{\sqrt{\pi}}{\sqrt{2}\,3^{1/4}}r+\tfrac{3^{5/6}\pi^{1/3}}{2^{4/3}}r^{2/3}+O(r^{1/3}).
\]
\end{theorem}
\begin{IEEEproof}
The proof is the same as the proof of
Theorem~\ref{thm:Euc_upper_square}, using
Lemmas~\ref{lem:Euc_circ_hex} and~\ref{lem:Euc_anti_hex} in place of
Lemmas~\ref{lem:Euc_circ_square} and~\ref{lem:Euc_anti_square}, and
defining $\ell=\lfloor 3^{1/12}2^{-5/6}\pi^{-1/6}r^{2/3}\rfloor$.
\end{IEEEproof}

\subsection{More general shapes}

All the theorems above consider a maximal anticode in some metric,
and cover this region with small circles of radius $\ell$. We comment
(for use later) that the same techniques work for any `sensible' shape
that is not necessarily an anticode. (We just need that the number of
small circles that intersect our shape is approximately equal to the
number of grid points contained in the shape.) So we can prove similar
theorems for DDCs that are restricted to lie inside regular polygons,
for example. The maximal number of dots in such a DDC is at most
$\sqrt{s} +o(\sqrt{s})$ when the shape contains $s$ points of the
grid.

\section{Periodic Two-Dimensional Configurations}
\label{sec:TwoDPatterns}

The previously known constructions for DDCs restrict dots to lie in a
line or a rectangular region (often a square region) of the plane.  The
application described in \cite{BEMP} instead demands that the dots lie in some anticode. The
most straightforward approach to constructing DDCs for our application
is to find a large square or rectangular subregion of our anticode,
and use one of these known constructions to place dots in this
subregion. This approach provides a lower bound for $m$ that is linear
in $r$, but in fact we are able to do much better than this by
modifying known constructions (in the case of Robinson's folding
technique below) and by making use of certain periodicity properties
of infinite arrays related to rectangular constructions. We will
explain how this can be done in the next section. In this section we
will survey some of the known constructions for rectangular DDCs,
extend these constructions to infinite periodic arrays, and prove the
properties we need for Section~\ref{sec:lower}.

Let $\cA$ be a (generally infinite) array of dots in the square grid,
and let $\eta$ and $\kappa$ be positive integers. We say that $\cA$ is
\emph{doubly periodic} with period $(\eta,\kappa)$ if
$\cA(i,j)=\cA(i+\eta,j)$ and $\cA(i,j)=\cA(i,j+\kappa)$ for all
integers $i$ and $j$. We define the \emph{density} of $\cA$ to be
$d/(\eta\kappa)$, where $d$ is the number of dots in any $\kappa
\times \eta$ sub-array of $\cA$. Note that the period $(\eta,\kappa)$
will not be unique, but that the density of $\cA$ does not depend on
the period we choose. We say that a doubly periodic array $\cA$ of
dots is \emph{a doubly periodic $n\times k$ DDC} if every $n\times k$
sub-array of $\cA$ is a DDC. See~\cite{Etz89,MGT93,Tay84} for some
information on doubly periodic arrays in this context. We aim to
present several constructions of doubly periodic DDCs of high density.

\subsection{Constructions from Costas Arrays}

A  {\it Costas array} of order $n$ is an $n \times n$ permutation
array which is also a DDC. Essentially two constructions for Costas
arrays are known, and both give rise to doubly periodic DDCs.

\vspace{3mm}
\noindent {\bf The Periodic Welch Construction:}

Let $\alpha$ be a primitive root modulo a prime $p$ and let $\cA$
be the square grid. For any integers $i$ and $j$, there is a dot
in $\cA(i,j)$ if and only if $\alpha^i \equiv j \bmod p$.

The following theorem is easy to prove. A proof which also mentions
some other properties of the construction is given
in~\cite{BEMP1}.

\begin{theorem}
Let $\cA$ be the array of dots from the Periodic Welch
Construction. Then $\cA$ is a doubly periodic $p\times (p-1)$ DDC with
period $(p-1,p)$ and density $1/p$.
\end{theorem}
Indeed, it is not difficult to show that each $p \times (p-1)$
sub-array is a DDC with $p-1$ dots: a dot in each
column and exactly one empty row. The $(p-1) \times (p-1)$ sub-array with
lower left corner at $\cA(1,1)$ is a Costas array.

\vspace{3mm}
\noindent {\bf The Periodic Golomb Construction:}

Let $\alpha$ and $\beta$ be two primitive elements in GF($q$), where
$q$ is a prime power. For any integers $i$ and $j$, there is a dot in
$\cA(i,j)$ if and only if $\alpha^i + \beta^j =1$.

The following theorem is proved similarly to the proof
in~\cite{BEMP1,Gol84}.

\begin{theorem}
\label{thm:periodicGolomb}
Let $\cA$ be the array of dots from the Periodic Golomb Construction.
Then $\cA$ is a doubly periodic $(q-1)\times (q-1)$ DDC with period
$(q-1,q-1)$ and density $(q-2)/(q-1)^2$.
\end{theorem}
Indeed, each $(q-1) \times (q-1)$ sub-array of $\cA$ is a DDC with $q-2$ dots;
exactly one row and one column are empty. The $(q-2) \times (q-2)$ sub-array
with lower left corner at $\cA(1,1)$ is a Costas array.

If we take $\alpha = \beta$ in the Golomb construction, then the
construction is known as the Lempel Construction. There are various
variants for these two constructions resulting in Costas arrays with
orders slightly smaller (by 1, 2, 3, or 4) or larger by one than the
orders of these two constructions (see~\cite{Gol92,GoTa84}). These are
of less interest in our discussion, as they do not extend to doubly
periodic arrays in an obvious way.

%

\subsection{Constructions from Golomb rectangles}
\label{sec:Grectang}

A {\it Golomb rectangle} is an $n \times k$ DDC with $m$ dots; Costas
arrays are a special case. Apart from constructions of special cases, there is
essentially one other general construction known, the {\it folded
rulers} construction due to Robinson~\cite{Rob97}.

\noindent {\bf Folded Ruler Construction:}

Let $S= \{ a_1 , a_2 , \cdots , a_m \}\subseteq\{0,1,\ldots,n\}$ be a
Golomb ruler of length $n$. Let $\ell$ and $k$ be integers such that
$\ell \cdot k \leq n+1$. Define $\cA$ to be the $\ell \times k$ array
where $\cA(i,j)$, $0 \leq i \leq k-1$, $0 \leq j \leq \ell-1$, has a
dot if and only if $i \cdot \ell +j = a_t$ for some $t$.

\begin{theorem}
\label{thm:folding} The array $\cA$ of the Folded Ruler
Construction is an $\ell \times k$ Golomb rectangle.
\end{theorem}

We now show how to adapt the Folded Ruler Construction to obtain a
doubly periodic $\ell \times k$ DDC. We require a stronger object than
a Golomb ruler as the basis for our folding construction, defined as
follows.

\begin{definition}
Let $A$ be an abelian group, and let $D=\{a_1,a_2,\ldots,a_m\}\subseteq A$
be a sequence of $m$ distinct elements of $A$. We say that $D$ is a
\emph{$B_2$-sequence over $A$} if all the sums $a_{i_1}+a_{i_2}$ with
$1\leq i_1\leq i_2\leq m$ are distinct.
\end{definition}

For a survey on $B_2$-sequences and their generalizations the reader
is referred to~\cite{Bry04}. The following lemma is well known and can
be readily verified.

\begin{lemma}
\label{lem:B2_diff} A subset $D=\{a_1,a_2,\ldots,a_m\}\subseteq A$ is a
$B_2$-sequence over $A$ if and only if all the differences
$a_{i_1}-a_{i_2}$ with $1\leq i_1 \neq i_2\leq m$ are distinct in $A$.
\end{lemma}

So, in particular, a Golomb ruler is exactly a $B_2$-sequence over
$\bZ$. Note that a $B_2$-sequence $\{a_1,a_2,\ldots ,a_m\}$ over
$\bZ_n$ produces a Golomb ruler $\{b_1,b_2,\ldots,b_m\}$ whenever the
$b_i$ are integers such that $a_i\equiv b_i\bmod n$. Also note that if
$D$ is a $B_2$-sequence over $\bZ_n$ and $a\in\bZ_n$, then so is the
shift $a+D=\{a+d:d\in D\}$. The following theorem, due to
Bose~\cite{Bose42}, shows that large $B_2$-sequences over $\bZ_n$
exist for many values of $n$.

\begin{theorem}
\label{thm:Bose}
Let $q$ be a prime power. Then there exists a $B_2$-sequence
$a_1,a_2,\ldots,a_m$ over $\bZ_n$ where $n=q^2-1$ and $m=q$.
\end{theorem}

\noindent {\bf The Doubly Periodic Folding Construction:}

Let $n$ be a positive integer and $D= \{a_1 , a_2 , \ldots , a_m \}$
be a $B_2$-sequence in $\Z_n$. Let $\ell$ and $k$ be integers such
that $\ell \cdot k \leq n$. Let $\cA$ be the square grid.  For any
integers $i$ and $j$, there is a dot in $\cA(i,j)$ if and only if $a_t
\equiv i \cdot \ell + j \bmod n$ for some $t$.

\begin{theorem}
Let $\cA$ be the array of the Doubly Periodic Folding
Construction. Then $\cA$ is a doubly periodic $\ell \times k$ DDC of
period $(\frac{n}{\gcd(n,\ell)},n)$ and density $m/n$.
\end{theorem}
\begin{IEEEproof}
Let $f(x,y)=x \cdot \ell +y$. The period of $\cA$ follows from the
observation that for each two integers $\alpha$ and $\beta$ we have
$f(i,j)=f(i+\alpha \frac{n}{\gcd(n,\ell)} ,j + \beta n) \equiv i
\cdot \ell +j \bmod n$. The density of $\cA$ is $m/n$ follows since
there are exactly $m$ dots in any $n$ consecutive positions in any
column.

Let $\cS$ be an $\ell\times k$ sub-array, whose lower left-hand corner
is at $\cA(i,j)$. An alternative construction of the dots in $\cS$ is
as follows. Take the shift $(i \cdot \ell + j)+D$ of $D$, which is
also a $B_2$-sequence in $\bZ_n$. Let $D'$ be the corresponding Golomb
ruler in $\{0,1,\ldots ,n-1\}$, so $a\in D$ if and only if $a\equiv
b\bmod n$, where $b\in (i \cdot \ell + j)+D$. Then form dots in $\cS$
by using the Folded Ruler Construction. Hence, by
Theorem~\ref{thm:folding}, the dots in $\cS$ form a DDC and so the
theorem follows.
\end{IEEEproof}

The following slightly different construction also produces doubly
periodic Golomb rectangles.

\vspace{3mm}
\noindent {\bf The Chinese Remainder Theorem Construction:}

Let $n$ be a positive integer and let $D= \{a_1 , a_2 , \ldots , a_m
\}$ be a $B_2$-sequence in $\Z_n$. Let $n = \ell \cdot k$ be any
factorization of $n$ such that $\gcd(\ell,k)=1$. For any two
integers $i$ and $j$ we place a dot in $\cA(i,j)$, if and only if
$a_t= (i \cdot \ell + j \cdot k) \bmod n$ for some $t$.

\begin{theorem}
Let $\cA$ be the array constructed by the Chinese Remainder Theorem
construction. Then $\cA$ is a doubly periodic $\ell \times k$ DDC of
period $(k,\ell)$ and density $m/n$. Moreover, every $\ell\times k$
sub-array of $\cA$ contains exactly $m$ dots.
\end{theorem}
\begin{IEEEproof}
Let $f(x,y)=x \cdot \ell + y \cdot k$. For any two integers $\alpha$
and $\beta$ we have $f(i,j) \equiv f(i+\alpha k,j + \beta \ell) \equiv
i \cdot \ell + j \cdot k \bmod n$. So the definition of $\cA$ implies
that $\cA$ is doubly periodic with period $(k,\ell)$.

Since $\ell$ and $k$ are relatively primes, it follows (from the
Chinese Remainder Theorem) that each integer $s$ in the range $0 \leq
s \leq \ell \cdot k-1$, has a unique representation as $s=d \cdot \ell
+ e \cdot k$, where $0 \leq d \leq k-1$, $0 \leq e \leq \ell-1$. Hence
every $\ell \times k$ sub-array of $\cA$ has $m$ dots corresponding to
the $m$ elements of the $B_2$-sequence $D$. In particular, this
implies that $\cA$ has density $m/n$.

Assume for a contradiction that there exists an $\ell\times k$
sub-array $\cS$ of $\cA$ that is not a DDC. Suppose that the lower
left-hand corner of $\cS$ is at $\cA(i,j)$. The distribution of dots
in $\cS$ is the same as the distribution in the sub-array with lower
left-hand corner the origin once we replace $D$ by the shift
$(i\ell+j\ell)+D$. So, without loss of generality, we may assume that
the lower left-hand corner of $\cS$ lies at the origin. As the
distinct difference property fails to be satisfied, there are four
positions with dots in $\cA$ of the form:
\begin{displaymath}
\begin{array}{cc}
  \cA(i_1,j_1)~~~ & \cA(i_1+d,j_1+e) \\
  \cA(i_2,j_2)~~~ & \cA(i_2+d,j_2+e) \\
\end{array}%
\end{displaymath}
where $i_1,i_1+d,i_2,i_2+d \in\{0,1,\ldots, k-1\}$ and
$j_1,j_1+e,j_2,j_2+e \in\{0,1,\ldots ,\ell-1\}$. By the definition of
$\cA$ we have
\begin{displaymath}
\begin{array}{c}
(i_1+d) \ell+(j_1+e)k-(i_1 \ell+j_1 k)= d \cdot \ell + e \cdot k\\
(i_2+d) \ell+(j_2+e)k-(i_2 \ell+j_2 k)=d \cdot \ell + e \cdot k
\end{array}%
\end{displaymath}
Since by Lemma~\ref{lem:B2_diff} each nonzero residue $s$ modulo $n$
has at most one representation as a difference from two elements of
$D$, it follows that the pairs
\begin{displaymath}
\begin{array}{c}
 \{  \cA(i_1,j_1), \cA(i_1+d,j_1+e) \} \\
  \{ \cA(i_2,j_2) , \cA(i_2+d,j_2+e) \} \\
\end{array}%
\end{displaymath}
are identical, and the theorem follows.
\end{IEEEproof}


\section{Lower bounds}
\label{sec:lower}

\subsection{Manhattan distance}
\label{sec:optimal_Lee}

In this section we will prove that there exists a $\oDD(m,r)$ with
$\frac{r}{\sqrt{2}} -o(r)$ dots: this attains asymptotically the
upper bound of Theorem~\ref{thm:up_Lee}. We will see that this
construction is actually using folding in a slightly different way.
We further show that we can construct a doubly periodic array in
which each Lee sphere of diameter $r$ is a DDC with
$\frac{r}{\sqrt{2}} +o(r)$ dots.

\vspace{3mm}
\noindent {\bf The LeeDD Construction:}

Let $r$ be an integer, and define $R=\lfloor \frac{r}{2} \rfloor$. Let
$D=\{a_1,a_2,\ldots,a_\mu\}$ be a ruler of length $n$. Define
$f(i,j)=iR+j(R+1)+R^2+R$. Let $\cA$ be the Lee sphere of radius $R$
centred at $(0,0)$, so $\cA$ has the entry $\cA(i,j)$ if $|i|+|j|
\leq R$. We place a dot in $\cA(i,j)$ if and only if $f(i,j) \in D$.

\begin{theorem}
\label{thm:LeeDD} The Lee sphere $\cA$ of the LeeDD Construction
is a $\oDD(m,r)$, where $m=\left \lvert D\cap\{0,1,\ldots,2R^2+2R\}\right\rvert$.
\end{theorem}
\begin{IEEEproof}
We first note that if $|i|+|j| \leq R$ then the smallest value that
the function $f$ takes is 0 and the largest value is $2R^2+2R$.  Next,
we claim that if $(i_1,j_1)$ and $(i_2,j_2)$ are two distinct points
such that $|i_1|+|j_1| \leq |i_2|+|j_2| \leq R$ then $f(i_1,j_1) \neq
f(i_2,j_2)$. Assume the contrary, that $f(i_1,j_1) = f(i_2,j_2)$. So
$i_1R+j_1(R+1)+R^2+R=i_2R+j_2(R+1)+R^2+R$ and therefore
$(i_2-i_1)R=(j_1-j_2)(R+1)$. If $i_1=i_2$, then $j_1=j_2$ which
contradicts our assumption that $(i_1,j_1)$ and $(i_2,j_2)$ are
distinct. So we may assume that $i_1\not=i_2$. Similarly, we may
assume that $j_1\not=j_2$. The equality $(i_2-i_1)R=(j_1-j_2)(R+1)$
now implies that $R+1$ divides $|i_2-i_1|$ and $R$ divides
$|j_2-j_1|$. This implies that $|i_2-i_1|+|j_2-j_1| > 2R$,
but
\[
|i_2-i_1|+|j_2-j_1|\leq |i_1|+|j_1|+ |i_2|+|j_2| \leq 2R,
\]
and so we have a contradiction. Thus, $f(i_1,j_1) \neq
f(i_2,j_2)$. This implies that each one of the integers between $0$
and $2R^2+2R$ is the image of exactly one pair $(i,j)$. In particular,
the number $m$ of dots in the configuration is exactly
$\left\lvert D\cap\{0,1,\ldots,2R^2+2R\}\right\rvert$.

Since $\cA$ is a Lee sphere of radius $R$, it follows that the
Manhattan distance between any two points is at most $2R \leq r$.
Now, assume for a contradiction that $\cA$ is not a $\oDD(m,r)$, so there
exist four positions with dots in $\cA$ as follows:
\begin{displaymath}
\begin{array}{cc}
  \cA(i_1,j_1)~~~ & \cA(i_1+d,j_1+e) \\
  \cA(i_2,j_2)~~~ & \cA(i_2+d,j_2+e) \\
\end{array}%
\end{displaymath}
By definition we have that
\[
f(i_1,j_1), f(i_1+d,j_1+e),
f(i_2,j_2),f(i_2+d,j_2+e) \in D.
\]
But then
$f(i_1+d,j_1+e)-f(i_1,j_1)=f(i_1+d,j_1+e)-f(i_1,j_1)=dR+e(R+1)$,
contradicting the fact that $D$ is a ruler.

Thus, the Lee sphere $\cA$ of the LeeDD Construction is a
$\oDD(m,r)$.
\end{IEEEproof}
\begin{cor}
There exists a $\oDD(m,r)$ in which $m = \frac{r}{\sqrt{2}}-o(r)$.
\end{cor}
\begin{IEEEproof}
Define $R=\lfloor r/2\rfloor$ and let $n=2R^2+2R+1$. There exists a
ruler of length at most $n$ containing $m$ dots, where $m\geq
\sqrt{n}+o(\sqrt{n})$: see~\cite{ASU,LaSa88,RoBe67}. Let
$D\subseteq\{0,1,\ldots,n-1\}$ be such a ruler. The corollary now
follows, by Theorem~\ref{thm:LeeDD}.
\end{IEEEproof}

It is worth mentioning that the LeeDD Construction is actually a
folding of the ruler by the diagonals of the Lee
sphere. Figure~\ref{fig:Leef} illustrates why this is the case, by
labelling the positions in a Lee sphere of radius $3$ by the values of
$f(i,j)$ at these positions. So if we use a $B_2$-sequence over
$\bZ_n$ instead of a ruler in the LeeDD Construction we obtain a
doubly periodic array with nice properties:

\begin{figure}[t]
\centering
\setlength{\unitlength}{.85mm}
\begin{picture}(45,40)(0,5)
\put(5,20){\framebox(35,5){}}
\put(10,15){\framebox(25,15){}}
\put(15,10){\framebox(15,25){}}
\put(20,5){\framebox(5,35){}}

\put(20,35){\makebox(5,5){24}}

\put(15,30){\makebox(5,5){17}}
\put(20,30){\makebox(5,5){20}}
\put(25,30){\makebox(5,5){23}}

\put(10,25){\makebox(5,5){10}}
\put(15,25){\makebox(5,5){13}}
\put(20,25){\makebox(5,5){16}}
\put(25,25){\makebox(5,5){19}}
\put(30,25){\makebox(5,5){22}}

\put(5,20){\makebox(5,5){3}}
\put(10,20){\makebox(5,5){6}}
\put(15,20){\makebox(5,5){9}}
\put(20,20){\makebox(5,5){12}}
\put(25,20){\makebox(5,5){15}}
\put(30,20){\makebox(5,5){18}}
\put(35,20){\makebox(5,5){21}}

\put(10,15){\makebox(5,5){2}}
\put(15,15){\makebox(5,5){5}}
\put(20,15){\makebox(5,5){8}}
\put(25,15){\makebox(5,5){11}}
\put(30,15){\makebox(5,5){14}}

\put(15,10){\makebox(5,5){1}}
\put(20,10){\makebox(5,5){4}}
\put(25,10){\makebox(5,5){7}}

\put(20,5){\makebox(5,5){0}}
\end{picture}

\caption{Folding along diagonals}
\label{fig:Leef}
\end{figure}
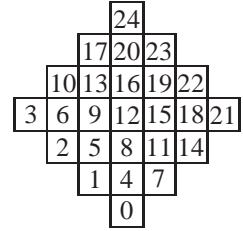

\vspace{3mm}
\noindent {\bf The Doubly Periodic LeeDD Construction:}

Let $r$ be an integer, $R=\lfloor \frac{r}{2} \rfloor$, and let
$D=\{a_1,a_2,\ldots,a_\mu\}$ be a $B_2$-sequence over $\Z_n$, where $n
\geq 2R^2+2R+1$. Let $f(i,j) \equiv iR+j(R+1) \bmod n$. Let $\cA$
be the square grid. For each two integers $i$ and $j$, there is a
dot in $\cA(i,j)$ if and only if $f(i,j) \in D$.

Similarly to Theorem~\ref{thm:LeeDD} we can prove the following
result.

\begin{theorem}
\label{thm:LeeDDthm}
The array $\cA$ constructed in the LeeDD Construction is doubly
periodic with period $(n,n)$ and density $\mu/n$. The dots contained
in any Lee sphere of radius $R$ form a DDC.
\end{theorem}
\begin{proof}
The first statement of the theorem is obvious. The second statement
follows as in the proof of Theorem~\ref{thm:LeeDD}, once we observe
that $f$ is an injection when restricted to any Lee sphere of radius
$R$.
\end{proof}

\begin{figure}
\centering
\setlength{\unitlength}{.6mm}
\begin{picture}(55,55)(0,0)
\multiput(0,0)(5,5){8}{
\multiput(0,0)(-5,5){4}{
\put(15,0){\line(1,0){5}}
\put(15,0){\line(0,1){5}}
\put(20,0){\line(0,1){5}}
\put(20,5){\line(-1,0){5}}
}
}
\end{picture}

\caption{A (3,5)-diagonally extended Lee sphere}
\label{fig:diagrec}
\end{figure}
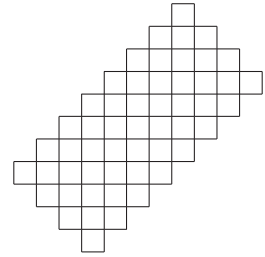

In Subsections~\ref{sec:hexd} and \ref{sec:euchex} we will make use of
an extension of this construction.  For positive integers $R$ and $t$,
an {\em $(R,t)$-diagonally extended Lee sphere} is a set
of positions in the square grid defined as follows. Let $(i_0,j_0)\in \bZ^2$,
and define $C=\{(i_0+k,j_0+k):0\leq k\leq t-1\}$. Then an
$(R,t)$-diagonally extended Lee sphere is the union of the Lee spheres
of radius $R$ with centres lying in $C$. (See Fig.~\ref{fig:diagrec}
for an example.) An $(R,t)$-diagonally extended Lee sphere contains
exactly $2R^2+t(2R+1)$ positions; the Lee sphere of radius $R$ is the
special case when $t=1$.
We observe that by choosing $n\geq 2R^2+t(2R+1)$, we can generalize the
doubly periodic LeeDD construction by continuing folding along
the diagonals of the rectangle.  This yields the following corollary,
which will prove useful in the construction of configurations for the
hexagonal grid.
\begin{corollary}\label{cor:drec}
Let $a$ be positive, and let $n$ be an integer such that $n\geq
(2+2a)R^2+aR$.  Consider the array $\cA$ constructed using the
doubly periodic LeeDD Construction. Then $\cA$ is a doubly periodic
array with density $\mu/n$. The dots contained in any $(R,\lfloor
aR\rfloor)$-diagonally extended Lee sphere form a DDC. There exists a
family of $B_2$ sequences so that $\cA$ has density at least $1/\sqrt{(2+2a)R^2+o(R^2)}$.
\end{corollary}
\begin{proof}
To establish the final statement of the corollary, we choose a family of
$B_2$ sequences as follows. Let $p$ be the smallest prime such that
$p^2-1\geq (2+2a)R^2+aR$, and define $n=p^2-1$. By Ingham's classical
result~\cite{Ingham37} on the gaps between primes, we have that
$n\leq (2+2a)R^2+O(R^{13/8})=(2+2a)R^2+o(R^2)$. By Theorem~\ref{thm:Bose}, there exists a $B_2$ sequence over $\bZ_n$ with $\mu=p$. Hence the density of $\cal A$ is
\[
\mu/n=p/(p^2-1)\geq1/\sqrt{(2+2a)R^2+o(R^2)},
\]
as required.
\end{proof}

\subsection{A General Technique}

Let $\cS$ be a shape (a set of positions) in the square grid. We are
interested in finding large DDCs contained in $\cS$, where (for
example) $\cS$ is an anticode. This subsection presents a general
technique for showing the existence of such DDCs, using the doubly
periodic constructions from Section~\ref{sec:TwoDPatterns}.

We write $(i,j)+\cS$ for the shifted copy
$\{(i+i',j+j'):(i',j')\in\cS\}$ of $\cS$. Let $\cA$ be a doubly
periodic array. We say that $\cA$ is a \emph{doubly periodic
$\cS$-DDC} if the dots contained in every shift $(i,j)+\cS$ of $\cS$
form a DDC. So the doubly periodic arrays constructed in
Section~\ref{sec:TwoDPatterns} are all doubly periodic $\cS$-DDCs
where $\cS$ is a square or a rectangle; the arrays in
Theorem~\ref{thm:LeeDDthm} and Corollary~\ref{cor:drec} are doubly
periodic $\cS$-DDCs with $\cal S$ a Lee sphere and diagonally extended
Lee sphere respectively. The following lemma follows in a
straightforward way from our definitions:

\begin{lemma}
\label{lem:subshape}
Let $\cA$ be a doubly periodic $\cS$-DDC, and let $\cS'\subseteq
\cS$. Then $\cA$ is a doubly periodic $\cS'$-DDC.
\end{lemma}

We will use doubly periodic DDCs to prove the existence of the
configurations we are most interested in, using the following theorem.

\begin{theorem}
\label{thm:general_construction}
Let $\cS$ be a shape, and let $\cA$ be a doubly periodic $\cS$-DDC of
density $\delta$. Then there exists a set of at least
$\left\lceil\delta\lvert\cS\rvert\right\rceil$ dots contained in $\cS$ that form a DDC.
\end{theorem}
\begin{IEEEproof}
Let the period of $\cA$ be $(\eta,\kappa)$. Write $m_{i,j}$ for the
number of dots of $\cA$ contained in the shift $(i,j)+\cS$ of
$\cS$. Now $\cA$ is periodic, so the definition of the density of
$\cA$ shows that
\[
\sum_{i=1}^\eta\sum_{j=1}^\kappa m_{i,j} = (\eta\kappa)\delta|\cS|.
\]
Hence the average size of the integer $m_{i,j}$ is $\delta\lvert\cS\rvert$, so
there exists an integer $m_{i',j'}$ such that $m_{i',j'}\geq
\left\lceil\delta\lvert\cS\rvert\right\rceil$. The $m_{i',j'}$ dots in $(i',j')+\cS$ form a
DDC, by our assumption on $\cA$, and so the appropriate shift of these
dots provides a DDC in $\cS$ with at least $\left\lceil\delta\lvert\cS\rvert\right\rceil$
dots, as required.
\end{IEEEproof}

\subsection{Euclidean distance in the square model}
\label{euclidsquare}
This subsection illustrates our general technique in the square grid
using Euclidean distance. So we wish to construct a $\DD(m,r)$ with
$m$ as large as possible.

\begin{figure}[tb]

\centering

\setlength{\unitlength}{.8mm}

\begin{picture}(80,70)(0,6)


\linethickness{.8 pt}
\put(40,40){\line(-100,168){15.359}}
\put(40,40){\line(100,168){15.359}}
\put(40,40){\line(0,1){25.7696}}
\put(29,53){\makebox(0,0){$R$}}
\put(38,48){\makebox(0,0){$\theta$}}

\put(40,42){\makebox(4,25){$\left.\rule{0cm}{12\unitlength}\right\}$}}
\put(45,55){\makebox(0,0){$\frac{n}{2}$}}

\put(40,40){\arc[90,120.76]{10}}


\put(24.6410,65.7679){\line(1,0){30.718}}

\put(24.6410,14.2321){\line(1,0){30.718}}

\put(14.2321,24.6410){\line(0,1){30.718}}

\put(65.7679,24.6410){\line(0,1){30.718}}

%
%
%
%
%

\linethickness{1.2 pt}
\put(40,40){\circle{60}}
\linethickness{.2 pt}
\put(14.2321,14.2321){\framebox(51.5358,51.5358){}}
\end{picture}
\caption{Square intersecting a circle}

\label{fig:sqr_circle}
\end{figure}
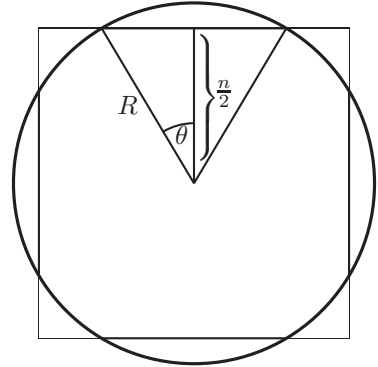

Let $R=\lfloor r/2\rceil$, and let $\cS$ be the set of points in the
square grid that are contained in the Euclidean circle of radius $R$
about the origin. We construct a DDC contained in $\cS$ with many
dots: any such configuration is clearly a $\DD(m,r)$ for some value of
$m$. The most straightforward approach is to find a large square
contained in $\cS$ (which will have sides of length approximately
$\sqrt{2}R$), and then add dots within this square using a Costas
array. This will produce a $\DD(m,r)$ where
\[
m=\sqrt{2}R-o(R)=\tfrac{1}{\sqrt{2}}r-o(r)\approx 0.707 r.
\]
To motivate our better construction, we proceed as follows. We find a
square of side $n$ where $n>\sqrt{2}R$ that partially overlaps our
circle: see Figure~\ref{fig:sqr_circle}. The constructions of
Section~\ref{sec:TwoDPatterns} show that there exist doubly periodic
$n\times n$ DDCs that have density approximately $1/n$. So
Theorem~\ref{thm:general_construction} shows that for any shape $\cS'$
within the square, there exist DDCs in $\cS'$ that have at least
$|\cS'|/n$ dots. Let $\cS'$ be the intersection of our square with
$\cS$. Defining $\theta$ as in the diagram, some basic geometry shows
that the area of $\cS'$ is
\[
|{\cal S'}| = \frac{(\pi/2)-2\theta+\sin 2\theta}{2\cos^2\theta}|{\cal S}|= 2R^2((\pi/2)-2\theta+\sin 2\theta).
\]
Since $n=2R\cos\theta$,
Theorem~\ref{thm:general_construction} shows that the density of dots
within $\cS'$ can be about $1/n=1/(2R\cos\theta)$ when $n$ is
large. So we can hope for at least $\mu R$ dots, where $\mu$ is the
maximum value of
\[
((\pi/2)-2\theta+\sin 2\theta)/\cos\theta
\]
on the interval $0\leq \theta\leq \pi/4$. In fact $\mu\approx 1.61589$,
achieved when $\theta\approx 0.41586$ (and so when $n=r\cos\theta=cr$,
where $c\approx 0.914769$).

\begin{theorem}
\label{thm:eucl_square_lower}
Let $\mu$ be defined as above. There exists a $\DD(m,r)$ in which
$m=(\mu/2)r-o(r)\approx 0.80795 r$.
\end{theorem}
Note that Theorem~\ref{thm:Euc_upper_square} gives an upper bound on
$m$ of the form $m\leq (\sqrt{\pi}/2)r+o(r)\approx 0.88623 r$.\newline
\begin{IEEEproof}
Define $c\approx 0.91477$ as above. Let $q$ be the smallest prime
power such that $q>cr$. We have that $cr<q<cr+(cr)^{5/8}$, by a
classical result of Ingham~\cite{Ingham37} on the gaps between primes;
so in particular $q\sim cr$. By Theorem~\ref{thm:periodicGolomb},
there exists a doubly periodic $(q-1)\times (q-1)$ DDC $\cA$ of density
$(q-2)/(q-1)^2$. Let $\cS'$ be the intersection between $\cS$ and a
Euclidean circle of radius $\lfloor r/2\rfloor$ about the origin. Then
$\cA$ is a doubly periodic $\cS'$-DDC. By
Theorem~\ref{thm:general_construction}, there exists a DDC in $\cS'$
with at least $m$ dots, where $|\cS'|(q-2)/(q-1)^2$. But the geometric
argument above shows that $|\cS'|(q-2)/(q-1)^2\sim (\mu/2)r$, and so the
theorem follows.
\end{IEEEproof}

\subsection{Hexagonal distance}
\label{sec:hexd}
By representing the hexagonal anticodes in the square grid, we may use
Theorem~\ref{thm:general_construction} to show the existence of a
$\oDD^*(m,r)$ where $m$ is large. The method of producing lower bounds
is essentially the same as above, but the geometrical problem we are
solving is different, with the images under $\xi$ of the maximal
anticodes $\cA_i$ replacing the circle, and the DDC contained a
diagonally extended Lee sphere of Corollary~\ref{cor:drec} replacing
the Costas array contained in a square.  Here we consider the case of
configurations contained in the hexagonal sphere
$\cA_{\lceil(r-1)/2\rceil}$; the cases of the other anticodes may be
handled in a similar fashion.
\begin{figure}[t]
\centering
\setlength{\unitlength}{1mm}
\begin{picture}(50,50)(-25,-25)
\linethickness{.2 pt}
\put(-20,-20){\dashbox(40,40){}}
\linethickness{.6 pt}
\put(11,20){\line(1,-1){9}}
\put(-20,-11){\line(1,-1){9}}

\put(-20,0){\line(-1,-1){5.5}}
\put(-20,-11){\line(-1,1){5.5}}
\put(20,0){\line(1,1){5.5}}
\put(20,11){\line(1,-1){5.5}}

\put(0,20){\line(1,1){5.5}}
\put(11,20){\line(-1,1){5.5}}
\put(0,-20){\line(-1,-1){5.5}}
\put(-11,-20){\line(1,-1){5.5}}

\put(0,20){$\underbrace{\rule{11\unitlength}{0cm}}$}
\put(2,11){\makebox(5,5){\footnotesize$aR$}}

\put(0,0){\line(-1,0){20}}
\put(-12,0){\makebox(5,5){\footnotesize$R$}}
\linethickness{1.2 pt}
\put(0,-20){\line(1,1){20}}
\put(-20,0){\line(1,1){20}}
\put(20,0){\line(0,1){20}}
\put(-20,0){\line(0,-1){20}}
\put(0,20){\line(1,0){20}}
\put(0,-20){\line(-1,0){20}}
\put(-10,-20){\makebox(7,7){\footnotesize $\cA_{\left\lceil\frac{r-1}{2}\right\rceil}$}}
\end{picture}
\caption{A diagonal rectangle intersecting the image of a hexagonal sphere}\label{fig:drechex}
\end{figure}
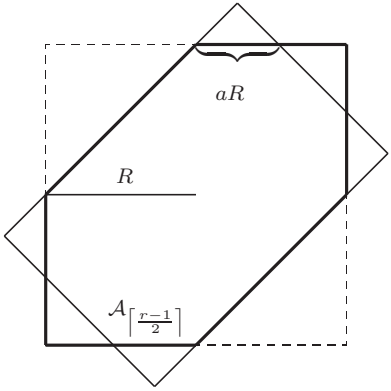
The problem we are solving is pictured in Fig.~\ref{fig:drechex}. The
figure shows the image under $\xi$ of the hexagonal sphere of radius
$R=\lfloor r/2\rfloor$ in bold; the square of side $2R+1$ containing this image is also
shown. The hexagonal sphere contains a Lee sphere of radius $R$ with
the same centre; the region $\cal S$ we consider is the $(R,\lfloor
aR\rfloor)$-diagonally extended Lee sphere whose mid-point is at the
centre of the hexagonal sphere: see Fig.~\ref{fig:drechex}. Let $\cal
S'$ be the intersection of $\cal S$ with the image of the hexagonal
sphere. We have that $|{\cal S'}|=R^2(2+2a-a^2)+o(R^2)$.  By
Corollary~\ref{cor:drec}, there is a doubly periodic $\cal S$-DDC of
density at least $1/\sqrt{n}$ where $n=2R^2(1+a)+o(R^2)$.  Thus
Theorem~\ref{thm:general_construction} shows that there is a DDC
contained in $\cal S'$ containing $\mu R-o(R)$ dots, where $\mu$ is the
maximum of
\begin{equation*}
\frac{2+2a-a^2}{\sqrt{2}\sqrt{1+a}}.
\end{equation*}
It can be seen that
$\mu=\left(\frac{2}{3}\right)^{\frac{3}{2}}\frac{1+2\sqrt{7}}{\sqrt{2+\sqrt{7}}}\approx1.58887$,
achieved when $a=\frac{-1+\sqrt{7}}{3}$.  Since $\cal S'$ is contained
in a hexagonal sphere of radius $R$, all pairs of dots in our DDC are at
hexagonal distance at most $r$. Thus we have the following theorem:
\begin{theorem}
Let $\mu$ be defined as above. There exists a $\overline{\DD}^*(m,r)$ in which
$m=(\mu/2)r-o(r)\approx 0.79444 r$.
\end{theorem}

\subsection{Euclidean distance in the hexagonal model}
\label{sec:euchex}
\begin{figure}
\centering
\setlength{\unitlength}{.5mm}
\begin{picture}(85,85)(-30,0)
\put(2.5,32.5){\vector(-1,0){30}}
\put(2.5,32.5){\vector(1,0){30}}
\put(-3,32.5){\makebox(10,5){\small $\mathbf 2t$}}

\put(45,45){\vector(1,1){7.5}}
\put(45,45){\vector(-1,-1){12.5}}
\put(50,37){\makebox(5,5){\small $\mathbf{\sqrt{2}at}$}}

\linethickness{.1pt}
\multiput(0,0)(5,5){11}{
\multiput(0,0)(-5,5){7}{
\put(0,0){\line(1,0){5}}
\put(5,0){\line(0,1){5}}
\put(0,0){\line(0,1){5}}
\put(0,5){\line(1,0){5}}
}
}

\end{picture}
\setlength{\unitlength}{.25mm}
\begin{picture}(155,155)(0,-60)

\put(65,2.886751345){\vector(-1,0){60}}
\put(65,2.886751345){\vector(1,0){60}}
\put(50,7.886751345){\makebox(40,10){\small $\mathbf 2t$}}
\put(50,-23.0947742){\vector(173,-100){45}}
\put(50,-23.0947742){\vector(-173,100){45}}
\put(110,-23.0947742){\vector(100,173){15}}
\put(110,-23.0947742){\vector(-100,-173){15}}
\put(40,-37){\makebox(5,5){\small $\mathbf{ \sqrt{3}t}$}}
\put(120,-23){\makebox(5,5){\small$\mathbf t$}}

\put(140,28.8675121){\vector(100,173){5}}
\put(140,28.8675121){\vector(-100,-173){15}}
\put(145,17){\makebox(5,5){\small$\mathbf{at}$}}

\linethickness{.1pt}

\multiput(0,0)(5,8.66025404){11}{
\put(0,0){\line(0,1){5.77350269}}
\put(0,5.77350269){\line(1732,1000){5}}
}

\multiput(0,0)(5,8.66025404){11}{
\multiput(0,0)(15,-8.66025404){7}{
\put(0,0){\line(1732,-1000){5}

\put(5,-2.88675135){\line(1732,1000){5}}
}}

\multiput(0,0)(15,-8.66025404){6}{
\put(15,-2.88675135){\line(-1732,1000){5}}
\put(15,-2.88675135){\line(0,-1){5.77350269}}
}}

\multiput(0,0)(5,8.66025404){10}{
\multiput(0,0)(15,-8.66025404){6}{
\put(10,0){\line(0,1){5.77350269}}
\put(15,-2.88675135){\line(1732,1000){5}}
}}

\multiput(95,-60.6217783)(5,8.66025404){11}{
\put(5,8.66025404){\line(0,1){5.77350269}}
\put(0,5.77350269){\line(1732,1000){5}}
}

\multiput(55,95.2627944)(15,-8.66025404){7}{
\put(0,0){\line(1732,-1000){5}}
}

\multiput(55,95.2627944)(15,-8.66025404){6}{
\put(5,-2.88675135){\line(0,-1){5.77350269}}
\put(15,-8.66025404){\line(-1732,-1000){5}}
}

\end{picture}
\caption{A $(t,\lfloor at\rfloor)$-diagonally extended Lee sphere is transformed into a rotated square (when $at=(\sqrt{3}-1)t+1$)}
\label{fig:rectangle_transformation}
\end{figure}
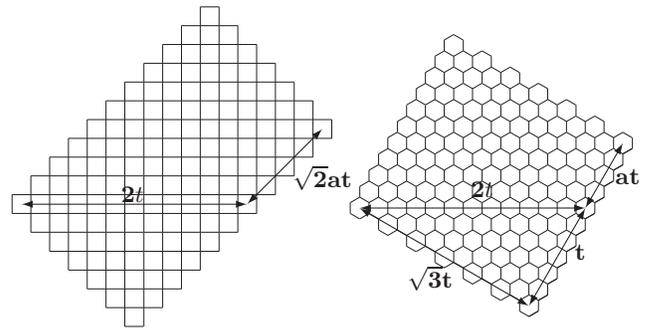

In this subsection we will obtain a construction for a ${\DD^*(m,r)}$
contained within a circle of radius $R=\lfloor r/2\rfloor$, again
based on the doubly periodic LeeDD construction.  We first observe
that a diagonally extended Lee sphere in the square grid is
transformed by $\xi^{-1}$ into a (rotated) rectangle in the hexagonal grid. In
particular, a $(t,\lfloor (\sqrt{3}-1)t+1 \rfloor)$-diagonally
extended Lee sphere is transformed by $\xi^{-1}$ into a set of
hexagons whose centres all lie within a (rotated) square $\cal S$ of side
$\sqrt{3}\,t$ (see
Fig.~\ref{fig:rectangle_transformation}). Corollary~\ref{cor:drec}
shows that there is a doubly periodic $\cal S$-DDC with density
$1/\sqrt{n}$, where $n=2\sqrt{3}t^2+o(t^2)$.

Consider (see Fig.~\ref{fig:sqr_circle2}) a circle of radius $R$ and a
square $\cal S$ of side $s$ where $s=2R\cos\theta$. Since a hexagon
has area $\sqrt{3}/2$, the square $\cal S$ contains
$(8/\sqrt{3})R^2\cos^2\theta + O(R)$ hexagons. Let $\cal S'$ be the
intersection of $\cal S$ with the circle of radius $R$. The
calculations in Subsection~\ref{euclidsquare} show that
\[
|{\cal S'}|=\frac{(\pi/2-2\theta+\sin 2\theta)}{2\cos^2\theta}|{\cal S}|+O(R).
\]
The previous paragraph shows that there is an periodic $\cal S'$-DCC
of density $\delta=1/\sqrt{n}$, where $n=(2/\sqrt{3})s^2+o(s^2)$. So
Theorem~\ref{thm:general_construction} now implies that there exists a
distinct difference configuration in $\cal S'$ containing at least $m$
dots, where
\[
m=\frac{\sqrt{\frac{2}{\sqrt{3}}}(\pi/2-2\theta+\sin
2\theta)}{\cos\theta}R-o(R).
\]
As in Subsection~\ref{euclidsquare}, we may take $\theta\approx
0.41586$ to maximise this expression. Hence we have proved the
following theorem:

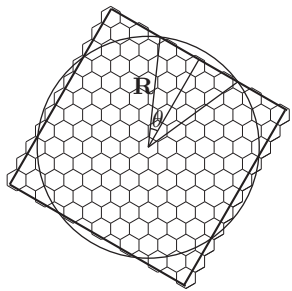
\begin{figure}[tb]
\centering
\setlength{\unitlength}{.25mm}
\begin{picture}(155,155)(0,-60)
\linethickness{.8pt}
\put(1.23463,.4422594){\line(1732,-1000){93.53}}
\put(1.23463,.4422594){\line(1000,1732){54}}
\put(55.23463,93.973){\line(1732,-1000){93.53}}
\put(93.7654,-53.5577){\line(1000,1732){54}}
\linethickness{.5pt}
\put(75,20.2072594){\line(1000,1732){27}}
\put(75,20.2072594){\line(100,924){6.4}}
\put(75,20.2072594){\line(1000,731){47.653}}

\put(75,20.2072594){\arc[60,83.827]{10}}
\put(77.5,33.2072594){\makebox(5,5){\small$\mathbf{\theta}$}}
\put(70,50.2072594){\makebox(5,5){\small$\mathbf{R}$}}

\put(75,20.2072594){\circle{118.06}}
\linethickness{.1pt}
\multiput(0,0)(5,8.66025404){11}{
\put(0,0){\line(0,1){5.77350269}}
\put(0,5.77350269){\line(1732,1000){5}}
}

\multiput(0,0)(5,8.66025404){11}{
\multiput(0,0)(15,-8.66025404){7}{
\put(0,0){\line(1732,-1000){5}

\put(5,-2.88675135){\line(1732,1000){5}}
}}

\multiput(0,0)(15,-8.66025404){6}{
\put(15,-2.88675135){\line(-1732,1000){5}}
\put(15,-2.88675135){\line(0,-1){5.77350269}}
}}

\multiput(0,0)(5,8.66025404){10}{
\multiput(0,0)(15,-8.66025404){6}{
\put(10,0){\line(0,1){5.77350269}}
\put(15,-2.88675135){\line(1732,1000){5}}
}}

\multiput(95,-60.6217783)(5,8.66025404){11}{
\put(5,8.66025404){\line(0,1){5.77350269}}
\put(0,5.77350269){\line(1732,1000){5}}
}

\multiput(55,95.2627944)(15,-8.66025404){7}{
\put(0,0){\line(1732,-1000){5}}
}

\multiput(55,95.2627944)(15,-8.66025404){6}{
\put(5,-2.88675135){\line(0,-1){5.77350269}}
\put(15,-8.66025404){\line(-1732,-1000){5}}
}

\end{picture}
\caption{Rotated square intersecting a circle}

\label{fig:sqr_circle2}
\end{figure}
\begin{theorem}
Let $\mu\approx 1.61589$ be the constant defined above
Theorem~\ref{thm:eucl_square_lower}.  There exists a $\rm{DD}^*(m,r)$
in which the number of dots is at least $\sqrt{\frac{2}{\sqrt{3}}}\mu
R-o(R)\approx 0.86819r$.
\end{theorem}

\section{Conclusion}
\label{sec:conclude}
We introduced the concept of a distinct difference configuration and
gave specific examples for both the square and hexagonal grids for
small parameters. We went on to provide general constructions for such
configurations, as well as upper and lower bounds on the maximum
number of dots such configurations may contain.  In the case of
distinct difference configurations using Manhattan distance these
bounds are tight asymptotically, as we have provided a construction
for configurations which meets the leading term in our upper bound.  For
the remaining classes of configurations, there is a gap between the
upper and lower bounds we have provided (see
Table~\ref{tab:boundsummary}). We believe the upper bounds to be
realistic, and it is an interesting challenge to provide constructions
that meet these bounds.  \renewcommand{\arraystretch}{1.4}
\begin{table}
\caption{Upper and lower bounds on the number of dots in a distinct difference configuration}\label{tab:boundsummary}
\begin{equation*}
\begin{array}{lcc}
\hline
& \text{lower bound} & \text{upper bound} \\
\hline
\overline{\rm{DD}}(m,r) & (1/\sqrt{2})r-o(r)& (1/\sqrt{2})r+O(r^{2/3})\\
\rm{DD}(m,r) &0.80795r-o(r)& 0.88623r+O(r^{2/3})\\
\overline{\rm{DD}}^*(m,r) & 0.79444r -o(r)&0.86603 r+O(r^{2/3}) \\
\rm{DD}^*(m,r) &0.86819r-o(r)&0.95231r+O(r^{2/3})\\
\hline
\end{array}
\end{equation*}
\end{table}
%
%

\begin{IEEEbiographynophoto}{Simon R. Blackburn}
received his BSc in Mathematics from the University of Bristol in
1989, and his DPhil in Mathematics from the University of Oxford in
1992. Since then he has worked at Royal Holloway, University of London
as a Research Assistant (1992-95), an Advanced Fellow (1995-2000), a
Reader in Mathematics (2000-2003) and a Professor in Pure Mathematics
(2004-). He was Head of the Mathematics Department from 2004 to
2007. His research interests include cryptography, group theory, and
combinatorics with applications to computer science.
\end{IEEEbiographynophoto}

\begin{IEEEbiographynophoto}{Tuvi Etzion}
(M'89-SM'99-F'04) was born in Tel Aviv, Israel, in 1956. He received
the B.A., M.Sc., and D.Sc. degrees from the Technion - Israel
Institute of Technology, Haifa, Israel, in 1980, 1982, and 1984,
respectively.

From 1984 he held a position in the department of Computer Science at
the Technion, where he has a Professor position. During the years
1986-1987 he was Visiting Research Professor with the Department of
Electrical Engineering - Systems at the University of Southern
California, Los Angeles. During the summers of 1990 and 1991 he was
visiting Bellcore in Morristown, New Jersey.  During the years
1994-1996 he was a Visiting Research Fellow in the Computer Science
Department at Royal Holloway, University of London. He also had
several visits to the Coordinated Science Laboratory at University of
Illinois in Urbana-Champaign during the years 1995-1998, two visits to
HP Bristol during the summers of 1996, 2000, several visits to the
department of Electrical Engineering, University of California at San
Diego during the years 2000-2009, and to the Mathematics department at
Royal Holloway, University of London during the years 2007-2009. His
research interests include applications of discrete mathematics to
problems in computer science and information theory, coding theory,
and combinatorial designs.

Dr Etzion was an Associate Editor for Coding Theory for
the IEEE Transactions on Information Theory from 2006 till 2009.
\end{IEEEbiographynophoto}

\begin{IEEEbiographynophoto}{Keith M. Martin} joined the Information
Security Group at Royal Holloway, University of London as a lecturer
in January 2000. He received his BSc (Hons) in Mathematics from the
University of Glasgow in 1988 and a PhD from Royal Holloway in
1991. Between 1992 and 1996 he held a Research Fellowship in the
Department of Pure Mathematics at the University of Adelaide,
investigating mathematical modeling of cryptographic key distribution
problems. In 1996 he joined the COSIC research group of the Katholieke
Universiteit Leuven in Belgium where he was primarily involved in an
EU ACTS project concerning security for third generation mobile
communications. He has also held visiting positions at the University
of Wollongong, University of Adelaide and Macquarie
University. Keith's current research interests include cryptography,
key management and wireless sensor network security.

Prof.\ Martin is an Associate Editor for Complexity and Cryptography
for IEEE Transactions on Information Theory.
\end{IEEEbiographynophoto}

\begin{IEEEbiographynophoto}{Maura B. Paterson}
received a BSc from the University of Adelaide in 2002 and a PhD from
Royal Holloway, University of London in 2005.  She has worked as a
research assistant in the Information Security Group at Royal
Holloway, and is currently at the Department of Economics, Mathematics
and Statistics at Birkbeck, University of London.  Her research
interests include applications of combinatorics in information
security.
\end{IEEEbiographynophoto}

\begin{thebibliography}{10}

\bibitem{Bab}
W. C. Babcock, ``Intermodulation interference in radio systems,''
\emph{Bull. Sys. Tech. Journal}, pp.\,63--73, June 1953.

\bibitem{Gol72}
S. W. Golomb, ``How to number a graph,'' in \emph{Graph Theory and
Computing}, Academic press, pp.\,23--37, 1972.

\bibitem{Shearer}
J. B. Shearer, ``Golomb rulers,''
\url{http://www.research.ibm.com/people/s/ shearer/grule.html}.

\bibitem{ASU}
M. D. Atkinson, N. Santoro, and J. Urrutia, ``Integer sets with
distinct sums and differences and carrier frequency assignments
for nonlinear repeaters'', \emph{IEEE Transactions on
Communications}, vol.\,COM-34, pp.\,614--617, 1986.

\bibitem{LaSa88}
A. W. Lam and D. V. Sarwate, ``On optimum time-hopping patterns'',
\emph{IEEE Transactions on Communications}, vol.\,COM-36,
pp.\,380--382, 1988.

\bibitem{GoTa82}
S. W. Golomb and H. Taylor, ``Two-dimensional synchronization
patterns for minimum ambiguity'', \emph{IEEE Trans.\ Inform.\
Theory}, vol.\,IT-28, pp.\,600--604, 1982.

\bibitem{RoBe67}
J. P. Robinson and A. J. Bernstein ``A class of binary recurrent
codes with limited error propagation'', \emph{IEEE Trans.\ Inform.\
Theory}, vol.\,IT-13, pp.\,106--113, 1967.

\bibitem{Rob85}
J. P. Robinson ``Golomb rectangles'', \emph{IEEE Trans.\ Inform.\
Theory}, vol.\,IT-31, pp.\,781--787, 1985.

\bibitem{Rob97}
J. P. Robinson ``Golomb rectangles as folded ruler'', \emph{IEEE
Trans.\ Inform.\ Theory}, vol.\,IT-43, pp.\,290--293, 1997.

\bibitem{Rob00}
J. P. Robinson ``Genetic search for Golomb arrays'', \emph{IEEE
Trans.\ Inform.\ Theory}, vol.\,IT-46, pp.\,1170--1173, 2000.

\bibitem{Cos75}
J. P. Costas, ``Medium constraints on sonar design and
performance,'' in \emph{EASCON Conv. Rec.}, pp.\,68A--68I, 1975.

\bibitem{GRT87}
R. Gagliardi, J. Robbins, and H. Taylor, ``Acquisition sequences in
PPM communications'', \emph{IEEE Trans.\ Inform.\ Theory},
vol.\,IT-33, pp.\,738--744, September 1987.

\bibitem{Gam87}
R. A. Games, ``An algebraic construction of sonar sequences using
M-sequences,'' \emph{SIAM J. Algebraic and Discrete Methods},
vol.\,8, pp.\,753--761, October 1987.

\bibitem{MGT93}
O. Moreno, R. A. Games, and H. Taylor, ``Sonar sequences from Costas
arrays and the best known sonar sequences with up to 100 symbols'',
\emph{IEEE Trans.\ Inform.\ Theory}, vol.\,IT-39, pp.\,1985--1987,
September 1993.

\bibitem{BlTi88}
A. Blokhuis and H. J. Tiersma, ``Bounds for the size of radar
arrays'', \emph{IEEE Trans.\ Inform.\ Theory}, vol.\,IT-34,
pp.\,164--167, January 1988.

\bibitem{HaZe97}
J. Hamkins and K. Zeger, ``Improved bounds on maximum size binary
radar arrays'', \emph{IEEE Trans.\ Inform.\ Theory}, vol.\,IT-43,
pp.\,997--1000, May 1997.

\bibitem{ZhTu94}
Z. Zhang and C. Tu, ``New bounds for the sizes of radar arrays'',
\emph{IEEE Trans.\ Inform.\ Theory}, vol.\,IT-40, pp.\,1672--1678,
September 1994.

\bibitem{EGRT92}
P. Erd\H{o}s, R. Graham, I. Z. Ruzsa, and H. Taylor, ``Bounds for
arrays of dots with distinct slopes or lengths'',
\emph{Combinatorica}, vol.\,12, pp.\,39--44, 1992.

\bibitem{PeTa00}
R. E. Peile and H. Taylor, ``Sets of points with pairwise distinct
slopes'', \emph{Computers and Mathematics}, vol.\,39, pp.\,109--115,
2000.

\bibitem{Zha93}
Z. Zhang, ``A note on arrays of dots with distinct slopes'',
\emph{Combinatorica}, vol.\,13, pp.\,127--128, 1993.

\bibitem{LeTh95}
H. Lefmann and T. Thiele, ``Point sets with distinct distances'',
\emph{Combinatorica}, vol.\,15, pp.\,379--408, 1995.

\bibitem{BEMP}
S. R. Blackburn, T. Etzion, K. M. Martin, and M. B. Paterson,
``Efficient key predistribution for grid-based wireless sensor
networks,'' \emph{Lecture Notes in Computer Science}, vol. 5155,
pp.\,54--69, August 2008.

\bibitem{BEMP1} S. R. Blackburn, T. Etzion, K. M. Martin, and
M. B. Paterson, ``Distinct difference configurations: multihop paths
and key predistribution in sensor networks'', preprint.

\bibitem{GoTa84}
S. W. Golomb and H. Taylor ``Constructions and properties of
Costas arrays'', \emph{Proceedings of the IEEE}, vol.\,72,
pp.\,1143--1163, 1984.

\bibitem{CoSl93}
J. H. Conway and N. J. A. Sloane, {\sl Sphere Packings, Lattices,
and Groups,} New York: Springer-Verlag, 1993.



\bibitem{AAK01}
R.\ Ahlswede, H.\,K.\ Aydinian, and L.\,H.\ Khachatrian, ``On
perfect codes and related concepts,''\! \emph{Designs, Codes
Crypto.}, vol. 22, pp. 221-237,~2001.

\bibitem{AhKh96}
R. Ahlswede and L. H. Khachatrian, ``The complete
nontrivial-intersection theorem for systems of finite sets,''
\emph{J.\ Combin.\ Theory, Series A}, vol.\,76, pp.\,121--138, 1996.

\bibitem{AhKh98}
R. Ahlswede and L. H. Khachatrian, ``The diametric theorem in
Hamming spaces--optimal anticodes,'' \emph{Adv. Appl. Math.},
vol.\,20, pp.\,429--449, 1998.%

\bibitem{Del73}
P. Delsarte, ``An algebraic approach to association schemes of
coding theory'', \emph{Philips J. Res.}, vol.\,10, pp.\,1--97,
1973.

\bibitem{ESV06}
T. Etzion, M. Schwartz, and A. Vardy, ``Optimal tristance anticodes
in certain graphs,'' \emph{J.\ Combin.\ Theory, Series A},
vol.\,113, pp.\,189--224, 2006.

\bibitem{MaZh95}
W. J. Martin and X. J. Zhu, ``Anticodes for the Grassman and
bilinear forms graphs,'' \emph{Designs, Codes, and Crypt.}, vol.\,6,
pp.\,73--79, 1995.

\bibitem{ScEt02}
M. Schwartz, and T. Etzion, ``Codes and anticodes in the Grassman
graph,'' \emph{J.\ Combin.\ Theory, Series A}, vol.\,97,
pp.\,27--42, 2002.

\bibitem{GoWe70}
S. W. Golomb and L. R. Welch, ``Perfect codes in the Lee metric and
the packing of polyominos'', \emph{SIAM J. Appl. Math.}, vol.\,18,
pp.\,302--317, 1970.

\bibitem{BBV98}
M. Blaum, J. Bruck, and A. Vardy, ``Interleaving schemes for
multidimensional cluster errors'', \emph{IEEE Trans.\ Inform.\
Theory}, vol.\,IT-44, pp.\,730--743, March 1998.

\bibitem{Etz03}
T. Etzion, ``Tilings with generalized Lee spheres,'' in
\emph{Mathematical Properties of Sequences and Other Combinatorial
Structures}, J. S. No, H. Y. Song, T. Helleseth, and P. V. Kumar,
editors, Kluwer Academic Publishers, pp.\,181--198, 2003.

\bibitem{Lit86} J.E. Littlewood (B. Bollob\'as, Ed),
\emph{Littlewood's miscellany}, Cambridge University Press, Cambridge,
1986.

\bibitem{Erd41} P. Erd\H{o}s and P. Tur\'an, ``On a problem of Sidon
in additive number theory and some related problems'', \emph{J. London
Math. Soc.}, vol.\,16, pp.\,212--215, 1941.

\bibitem{Etz89}
T. Etzion, ``Combinatorial designs derived from Costas arrays,''
in \emph{Sequences}, R. M. Capocelli editor, New York, NY:
Springer Verlag, pp.\,208--227, 1989.

\bibitem{Tay84}
H. Taylor ``Non-attacking rooks with distinct differences'',
Communication Sciences Institute, University of Southern California,
Tech. Report CSI-84-03-02, March 1984.

\bibitem{Gol84}
S. W. Golomb, ``Algebraic constructions for Costas arrays,''
\emph{J.\ Combin.\ Theory, Series A}, vol.\,37, pp.\,13--21, 1984.

\bibitem{Gol92}
S. W. Golomb, ``The $T_4$ and $G_4$ constructions for Costas
arrays'', \emph{IEEE Trans.\ Inform.\ Theory}, vol.\,IT-38,
pp.\,1404--1406, 1992.

\bibitem{Bry04}
K. O'Bryant, ``A complete annotated bibliography of work related to
Sidon sequences'', \emph{The Electronic Journal of Combinatorics},
DS11, pp.\,1--39, July 2004.

\bibitem{Bose42} R. C. Bose, ``An affine analogue of Singer's theorem'',
\emph{J. Indian Math.\ Soc.\ (N.S.)}, vol. 6, pp. 1-15, 1942.

\bibitem{Ingham37} A.E. Ingham, ``On the difference between
  consecutive primes'', \emph{Quart. J. Math. Oxford (O.S.)} vol. 8,
  pp.  255-266, 1937.
\end{thebibliography}
\end{document}